\documentclass[11pt]{article}

\usepackage{titlesec}
\titlelabel{\thetitle.\quad}
\usepackage{amsthm}
\usepackage{amsmath,amsthm,amsfonts,amssymb,mathrsfs,bm,graphicx,stmaryrd,dsfont}
\usepackage[colorlinks=true,linkcolor=blue]{hyperref}
\usepackage[usenames,dvipsnames]{color}
\usepackage{fullpage}
\usepackage[american]{babel}
\usepackage[varg]{pxfonts}
\usepackage{prettyref}
\usepackage{microtype}
\hypersetup{
   citecolor=blue
}
\usepackage{float}
\usepackage{amsopn}

\usepackage{algorithmic}
\usepackage[ruled, algosection]{algorithm2e}
\ifpdf
\hypersetup{
  pdftitle={The Lanczos Algorithm Under Few Iterations: Concentration and Location of the Output},
  pdfauthor={Jorge Garza-Vargas, Archit Kulkarni}
}
\fi

\usepackage{amsmath,amsfonts,amssymb,mathrsfs,bm,graphicx,stmaryrd,dsfont, todonotes}
\usepackage{mathtools}
\usepackage{enumerate}

\DeclarePairedDelimiter{\norm}{\lVert}{\rVert}
\newcommand{\ep}{\varepsilon}
\author{Jorge Garza-Vargas \\ UC Berkeley \\ jgarzavargas@berkeley.edu  \and  Archit Kulkarni \\ UC Berkeley \\ akulkarni@berkeley.edu }

\title{\Large\textbf{The Lanczos Algorithm Under Few Iterations: Concentration and Location of the Output}}

\date{\today}

\begin{document}

\maketitle
\allowdisplaybreaks
\newtheorem{theorem}{Theorem}[section]
\newtheorem{definition}[theorem]{Definition}
\newtheorem{conjecture}[theorem]{Conjecture}
\newtheorem{lemma}[theorem]{Lemma}
\newtheorem{corollary}[theorem]{Corollary}
\newtheorem{proposition}[theorem]{Proposition}
\newtheorem{observation}[theorem]{Observation}
\newtheorem{remark}[theorem]{Remark}
\newtheorem{result}[theorem]{Result}
\newtheorem{example}[theorem]{Example}
\newtheorem{question}[theorem]{Question}
\numberwithin{equation}{section}
\begin{abstract}
 We study the Lanczos algorithm where the initial vector is sampled uniformly from $\mathbb{S}^{n-1}$. Let $A$ be an $n \times n$ Hermitian matrix.  We show that when run for few iterations, the output of Lanczos on $A$ is almost deterministic.  More precisely, we show that for any $ \varepsilon \in (0, 1)$ there exists $c >0$ depending only on $\varepsilon$ and a certain global property of the spectrum of $A$ (in particular, not depending on $n$) such that when Lanczos is run for at most $c \log n$ iterations, the output Jacobi coefficients deviate from their medians by $t$ with probability at most $\exp(-n^\varepsilon t^2)$ for $t<\Vert A \Vert$. We directly obtain a similar result for the Ritz values and vectors.  Our techniques also yield asymptotic results:  Suppose one runs Lanczos on a sequence of Hermitian matrices $A_n \in M_n(\mathbb{C})$ whose spectral distributions converge in Kolmogorov distance with rate $O(n^{-\varepsilon})$  to a density $\mu$ for some $\varepsilon > 0$.  Then we show that for large enough $n$, and for $k=O(\sqrt{\log n})$, the Jacobi coefficients output after $k$ iterations concentrate around those for $\mu$. The asymptotic setting is relevant since Lanczos is often used to approximate the spectral density of an infinite-dimensional operator by way of the Jacobi coefficients;  our result provides some theoretical justification for this approach.  

     In a different direction, we show that Lanczos fails with high probability to identify outliers of the spectrum when run for at most $c' \log n$ iterations, where again $c'$ depends only on the same global property of the spectrum of $A$. Classical results imply that the bound $c' \log n$ is tight up to a constant factor. 
\end{abstract}

\section{Introduction}

Eigenvalue problems are ubiquitous in science and engineering. However, most applications require analyzing matrices whose large dimension makes it impractical to exactly compute any important feature of their spectrum. It is for this reason that iterative randomized algorithms have proliferated in numerical linear algebra \cite{saad2011numerical,trefethen1997numerical}. 

In this context, iterative randomized  algorithms provide an approximation of the spectrum of the matrix in question, where the accuracy of the approximation improves as the number of iterations increases. For any such algorithm, it is natural to ask the following questions: 

\begin{enumerate}
\item[(Q1)] How much does the random output vary?
    \item[(Q2)] How many iterations are necessary and sufficient to obtain a satisfactory approximation?  
\end{enumerate}

The present work, theoretical in nature, addresses the above questions for one of the most widely used algorithms for eigenvalue approximation, namely, the Lanczos algorithm. Throughout the paper we assume exact arithmetic. 

\subsection{The Lanczos algorithm}

Recall that when run for $k$ iterations, the Lanczos algorithm  outputs a $k\times k$ tridiagonal matrix called the \emph{Jacobi matrix}. The nontrivial entries of  the  the Jacobi matrix are the \emph{Jacobi coefficients}, which we denote by $\alpha_k$ and $\beta_k$, and its eigenvalues are the \emph{Ritz values}, which we denote by $r_k$.  Oftentimes, the Jacobi coefficients and the Ritz values provide important information about the spectrum of the matrix. In particular, when $k=n$, the Ritz values are exactly the eigenvalues of $A$, and hence the full spectrum is recovered. However, in practice it is usually too expensive to perform $\Theta(n)$ iterations. 
\bigskip

\noindent \textbf{Outlying eigenvalues and Ritz values.} The success of the Lanczos algorithm resides to some extent in its ability to find the \emph{outliers} of the spectrum of the matrix $A$ with very few iterations.  By outliers, we mean the eigenvalues distant from the region in which the majority of the spectrum accumulates (the \emph{bulk}).  Hence, the algorithm is of particular interest in most applications in science and engineering \cite{saad2011numerical}.  
\bigskip

\noindent \textbf{Bulk spectrum and Jacobi coefficients.} Lanczos-type methods can also be used to approximate the global spectral density of large matrices, also known as density of states; for a survey of techniques see \cite{lin2016approximating}.  
In applied mathematics, large matrices can arise as discretizations of infinite-dimensional operators such as the Laplacian or as finite-dimensional representations of an infinite-dimensional Hamiltonian.  Computing the eigenvalues and Jacobi coefficients of the finite-dimensional operator then yields information about the infinite-dimensional operator and the underlying continuous system.  For an example, see \cite{shao2018structure}, or Section 7 of \cite{van2001computing} for numerical experiments and bounds for the Lanczos algorithm applied to an explicit discretized Laplace operator. 

In the setting described above, the Jacobi coefficients contain all the information of the spectral density of the infinite-dimensional operator in question and even the first few coefficients are of use. To give an example, in \cite{haydock1980recursive} the Haydock method (as it is termed today) was introduced. This method exploits the fact the resolvent of an operator admits a continued fraction expansion where the coefficients are precisely the Jacobi coefficients, and hence knowing these quantities is fundamental to understanding the spectral density of the operator---see 3.2.2 \cite{lin2016approximating} for a summary of the Haydock method. 

Using a slightly  different perspective, note that from the $k\times k$ Jacobi matrix of an operator one can obtain the $[k-1, k]$ Pad\'e approximation of its resolvent \cite{van2006pade}. In particular, knowing the $k\times k$ Jacobi matrix is enough to compute the first $2k-1$ moments of the spectral density of the infinite-dimensional operator.     
\bigskip
 
In applications, sophisticated modifications of the Lanczos algorithm are used \cite{golub1977block, calvetti1994implicitly, li2016thick}. Since the goal of the present paper is to introduce proof techniques and theoretical tools that have not been exploited previously, we only deal with the simplest version of the Lanczos algorithm and do not strive to obtain optimal constants in our bounds and theorems when providing answers for questions (1) and (2).

\subsection{Question (1): Our contributions} \label{sec:contr1}

\begin{figure}[H]
    \centering

\begin{tabular}{c c c}
 \includegraphics[width=3.9cm, height=3.9cm]{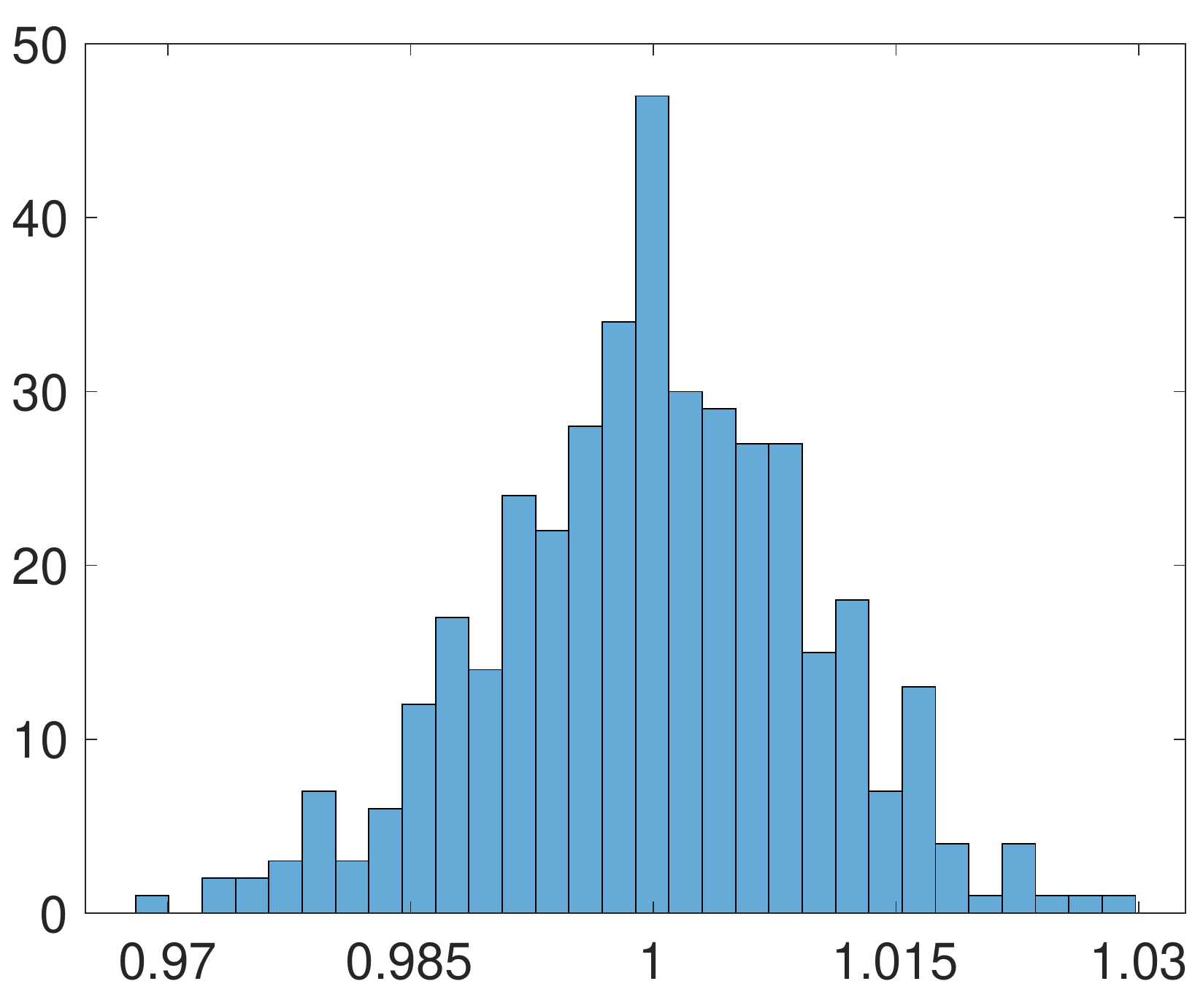}    & \includegraphics[width=3.85cm, height=3.85cm]{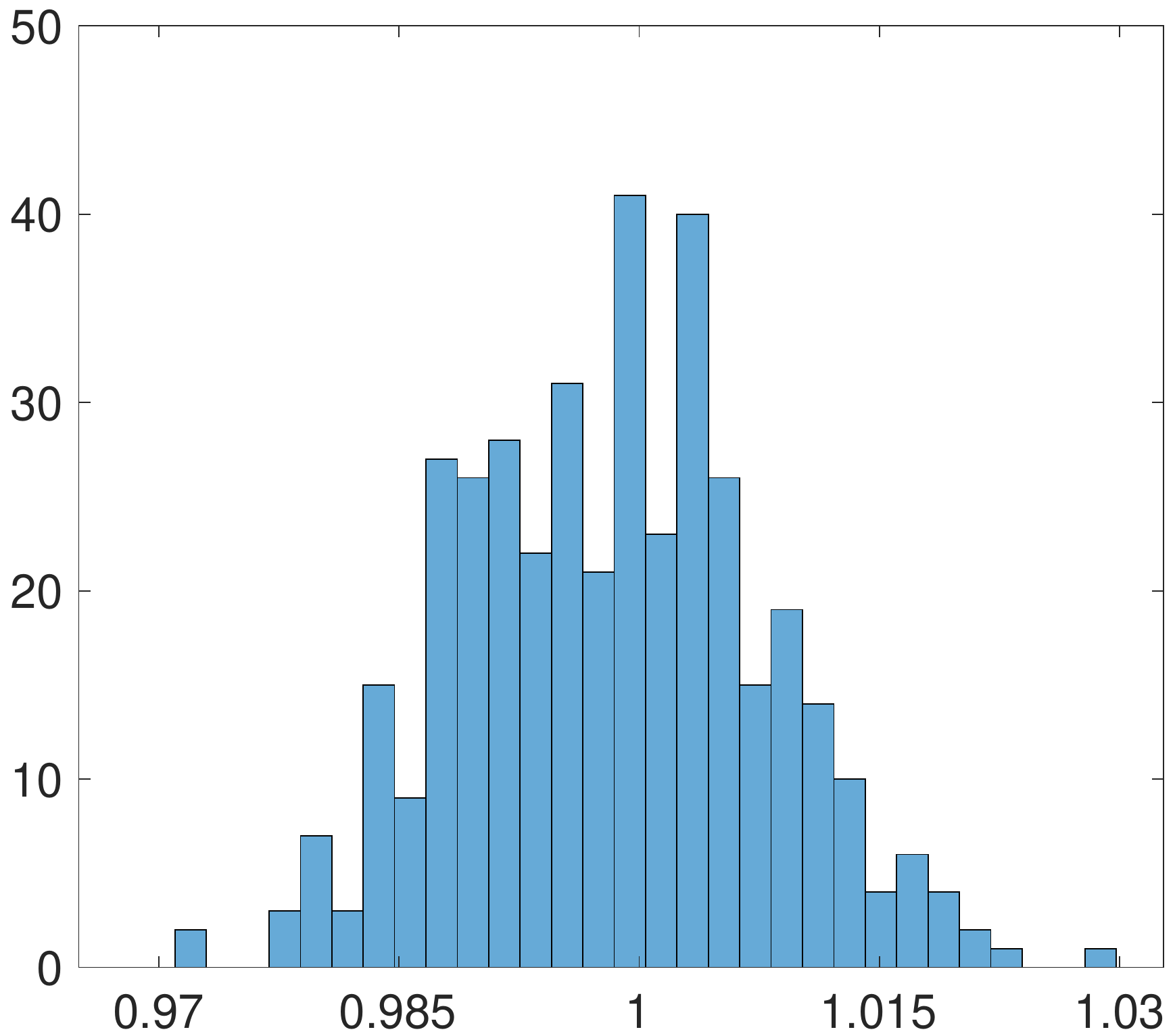}  & \includegraphics[width=3.85cm, height=3.85cm]{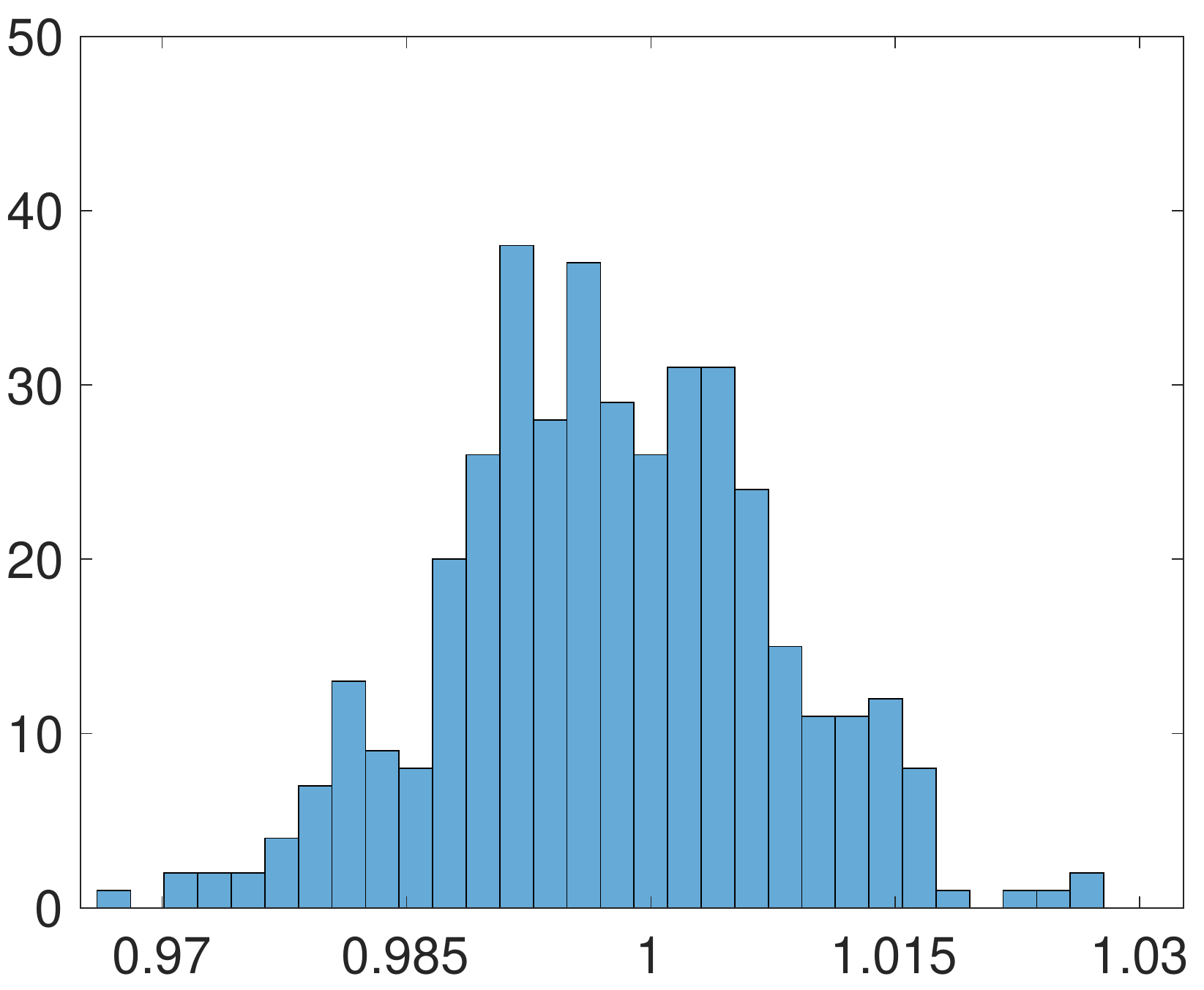} 
 \\ $\beta_0$ & $\beta_{10}$ & $\beta_{20}$
\end{tabular}
  \caption{ Here $A$ is a fixed $n \times n$ matrix drawn from the Gaussian orthogonal ensemble (GOE) with $n = 5000$. Since the empirical spectral distribution of $A$ will be close to the semicircular law it is expected that the Jacobi coefficients $\beta_i$ of $A$ will be approximately 1. The above histograms show the values $\beta_0, \beta_{10},$ and $\beta_{20}$ obtained by running  the Lanczos algorithm 400 times on the input $A$.   Note that in each of these cases, $\beta_i$ appears to be concentrated. }
    \label{fig:jacobicoefs}
\end{figure}

As far as we are aware, there is no previous work addressing this question for the Lanczos algorithm. In this paper we show that there is a $c >0$ such that for $n$ large enough, the output of the Lanczos procedure is almost deterministic when run for at most  $c \log n$ iterations. More precisely, in Theorem \ref{thmhessconc} we show that for $k\leq c\log(n)$ and $\varepsilon \in (0, 1/2)$,  deviations of the order $n^{-\varepsilon}$ of the Jacobi coefficients $\alpha_k$ and $\beta_k$ computed by Lanczos  occur with exponentially small probability. For an illustration, see Figure \ref{fig:jacobicoefs}. The strength of our probability bound deteriorates as $k$ grows.  The constants in the theorem depend only on an easily computed global property of the spectrum which we call \emph{equidistribution}.   

From the point of view of random matrix theory, the problem treated in the present paper is atypical.  In random matrix theory, most of the studied models have a rich probabilistic structure that can be exploited to obtain results about the eigenvalue distribution of the matrix. By contrast, in our case, the Jacobi matrix output by the Lanczos algorithm is a random matrix obtained by running a complicated deterministic dynamic over a minimal source of randomness---a single uniform random unit vector. Hence, in order to obtain results similar to the ones presented in this article, the structure of the algorithm needs to be exploited in an involved way. We use the ubiquitous concentration of measure phenomenon for Lipschitz functions in high dimension, together with a careful control of the variables appearing in the Lanczos algorithm and their Lipschitz constants as functions of the random input.  Roughly speaking, the Lipschitz constant is exponential in the number of iterations, which yields concentration in the regime of at most $c \log n$ iterations for sufficiently small $c$.  Throughout the analysis we use elementary results in the theory of orthogonal polynomials. 

In view of the fact that the output of the Lanczos algorithm is sharply concentrated under few iterations, one may ask which values the output is concentrated around.  Toward the end of this introduction we give an overview of our results in this direction.

\subsection{Question (2): Previous work}

Regarding the problem of detecting outliers of the spectrum via the Lanczos algorithm, theoretical answers to the sufficiency part of Question (2) posed above appeared decades ago. Most of them in essence give an \emph{upper bound} on the number of iterations required to obtain an accurate approximation when the input is an $n$-dimensional matrix $A$---see \cite{kaniel1966estimates, paige1971computation, saad1980rates}.  Roughly speaking, previous literature provides inequalities that state that $k \ge C \log n$ iterations suffice for the Lanczos algorithm to approximate the true extreme eigenvalues of $A$ very well, making the use of $O(\log n)$ iterations common in practice---see \cite{kuczynski1994probabilistic} or \cite{van2001computing} for examples of inequalities that give this bound. The constant $C$ in the results mentioned above is determined by features of the spectrum of $A$; typically, these features are the diameter of the spectrum and the gaps between the outliers and the bulk.  In recent years, more refined arguments have yielded inequalities in which other features of the spectrum are considered, see \cite{yuan2018superlinear} for an example or \cite{bellalij2010further} for a survey. 

Regarding the necessity part of Question (2), to the best of our knowledge, the only existing negative result regarding detection of outliers is the one given in the recent work \cite{simchowitz2018tight}. There,  a query complexity bound was proven for any algorithm that is allowed to make queries of matrix-vector products, which in particular applies to the Lanczos algorithm.

\subsection{Question (2): Our contributions}\label{sectionanswer2}

In the present paper we study the Lanczos algorithm in the context of approximation of outliers and answer the necessity part of Question (2).  That is, we show that if run for at most $k \le c \log n$ iterations,  the Lanczos algorithm fails to approximate outliers with overwhelming probability.  Thus, in essence we provide a \emph{lower bound} on the number of iterations required for accuracy.  As in our contribution for Question (1), the constant $c$ depends only on the equidistribution of the spectrum.
 
 To give some rough context, the result in \cite{simchowitz2018tight} discussed above shows that if the empirical spectral distribution of a matrix is close to the semicircle distribution plus an outlying ``spike,'' any algorithm in their class will fail to identify the spike with overwhelming probability, unless given at least $c \log n$ queries.  In contrast, our result applies exclusively to the Lanczos algorithm, but shows that outliers are missed for a far more general class of measures than just the semicircle.

 In order to analyze asymptotic behavior, we adopt a  framework similar to that used in 
 \cite{kuijlaars2000eigenvalues} and \cite{beckermann2000note}, in which a sequence of Hermitian matrices $A_n$ with convergent spectra was considered. These  papers studied the behavior of the Lanczos algorithm in the regime of $\Theta(n)$ iterations.

\begin{figure}[H]
\centering
\includegraphics[width=\textwidth]{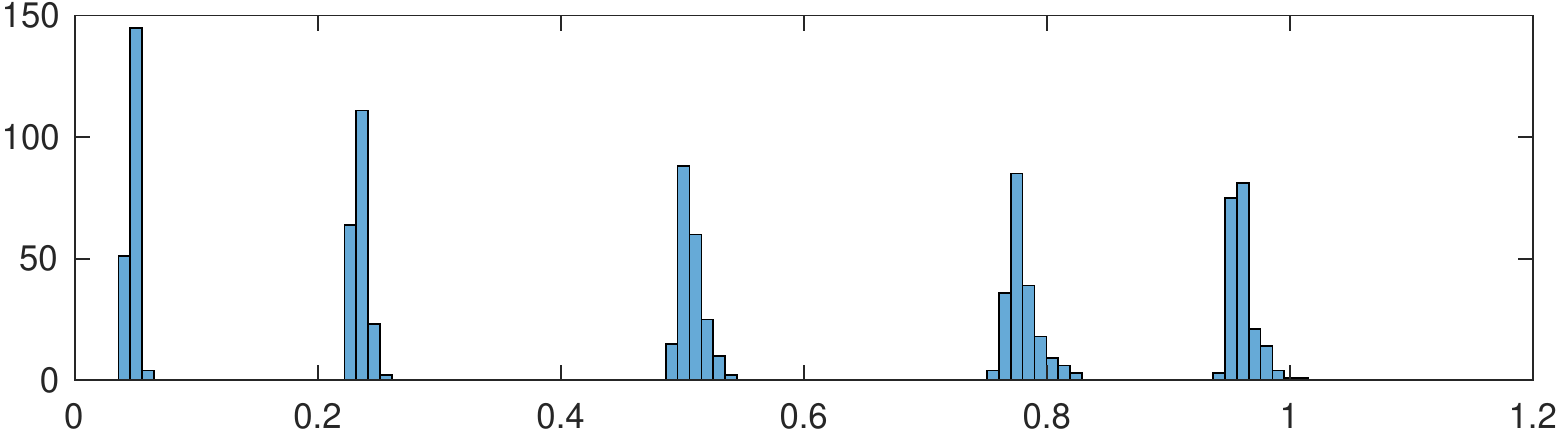}

\includegraphics[width=\textwidth]{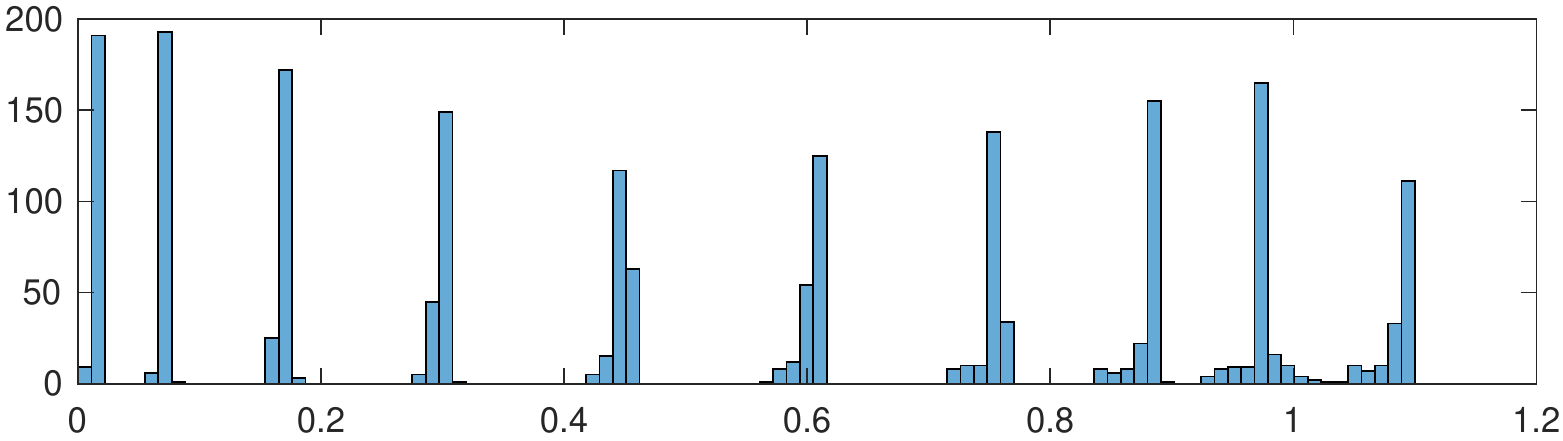}

\caption{$A$ is a $2000 \times 2000$ diagonal matrix with entries $\{0, 1/2000, 2/2000, \dots, 1999/2000, 1.1\}$. This represents a discretization of $\text{Unif}([0,1])$ plus an outlier at $1.1$.  Plotted is a histogram of the Ritz values output by Lanczos after $k=5$ iterations (above) and after $k=10$ iterations (below).  To generate the histogram the procedure was run 200 times. Notice that to find the outlier with a decent probability, 10 iterations suffice (but  5 do not).  However, even in the regime of $k=5$ iterations the output appears to be concentrated.}
\label{fig:outlier}
\end{figure}

To show that the Lanczos algorithm misses outliers when run for at most $c\log n$ iterations, we use the elementary theory of orthogonal polynomials  and standard techniques in high-dimensional probability. Roughly speaking, using a variational principle, we show that for small enough $k$, the roots of the $k$th orthogonal polynomial with respect to a certain random measure are contained in a small blow-up of the convex hull of the bulk of the true spectrum. See Theorem \ref{thm:nonasylocation} or Proposition \ref{prop:outliers}  for a precise statement and Figure \ref{fig:outlier} for an illustration. 

\subsection{Our result on the locations of the output} One may ask if finer statements about the location of the Jacobi coefficients and Ritz values can be made. Previously, tools from potential theory have been used to answer this question in the regime of $k=\Theta(n)$ iterations \cite{beckermann2000note, kuijlaars2000eigenvalues,kuijlaars2006convergence}. In the regime of $k$ fixed as $n \to \infty$, a deterministic convergence result for orthogonal polynomials \cite[Theorem 4]{gautschi1968construction} can be used to show that the Ritz values converge almost surely to the roots of the $k$th orthogonal polynomial of the limiting eigenvalue distribution; see Remark \ref{rem:fixedk} for details. 

\begin{figure}[H]

\centering

\includegraphics[width=\textwidth]{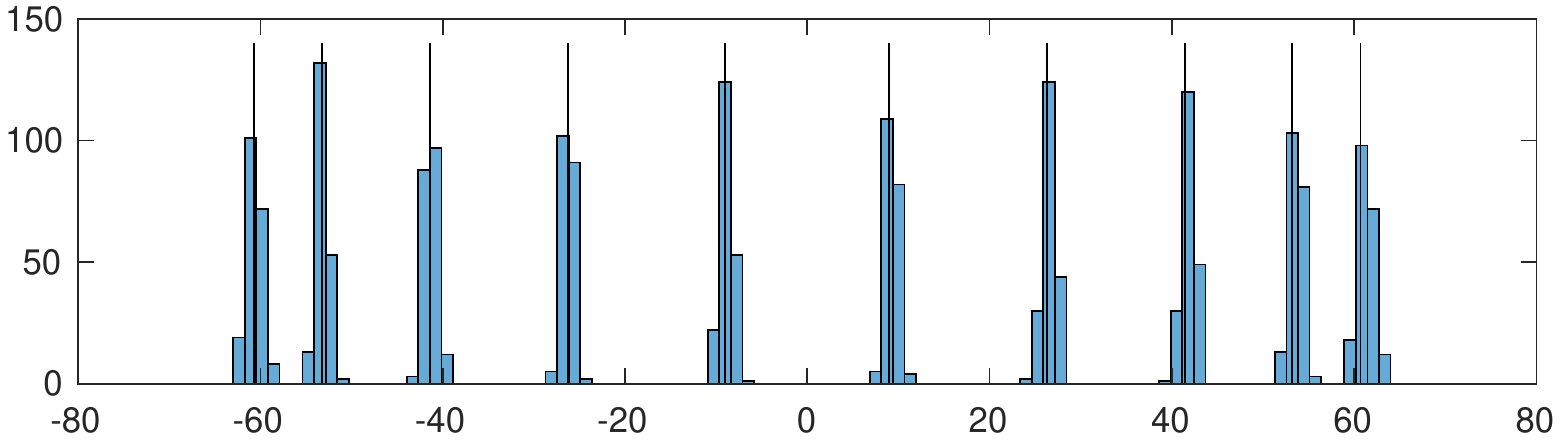}

\caption{$A$ is a fixed $n \times n$ matrix drawn from the GOE with $n = 2000$.  Plotted is the histogram of the Ritz values after 200 repetitions of the Lanczos algorithm with $k=10$ iterations.  Also plotted are the roots of the 10th orthogonal polynomial with respect to the (suitably rescaled) semicircle law, which is the limit of the distribution of eigenvalues for GOE as $n \to \infty$.}
\label{fig:wigner}
\end{figure}

In the present work we use determinantal formulas for orthogonal polynomials and concentration of measure results to locate the Jacobi coefficients and Ritz values in the regime of $k=O(\sqrt{\log n})$ iterations.  In particular, we prove that the Ritz values concentrate around the roots of the $k$th orthogonal polynomial for the limiting eigenvalue distribution.   See Figure \ref{fig:wigner} for an illustration. Moreover, also when $k= O(\sqrt{\log n})$, we show that the Jacobi matrix obtained after $k$  iterations is concentrated around the $k$th Jacobi matrix of the limiting measure. 

These results may be of particular relevance in applications where an infinite-dimensional operator is discretized with the goal of computing its density. In essence, Theorem \ref{Thmsqrtiter} below states that in this situation the first iterations of the Lanczos algorithm are an accurate approximation of the true Jacobi coefficients of the spectral measure of the infinite-dimensional operator, and hence the procedure is giving valuable information for recovering the limiting measure.

\subsection{Organization of the paper}

The article is organized as follows.  In Section 2, we review the classical background of the Lanczos procedure and orthogonal polynomials and formally state our main theorems.  In Section 3, we develop machinery that in Section 4 will be used to prove concentration for the output of the Lanczos algorithm.  In Section 5, we prove our complementary results about the location of the Ritz values and Jacobi coefficients. Finally, in Section 6 we discuss further research directions that may be of interest.

\section{Preliminaries and statements of theorems}

Throughout this paper only elementary facts about orthogonal polynomials are used. For the reader's convenience in Section 2.1 we include a concise survey of the results that will be used in what follows. Chapter 2 in \cite{szego1939orthogonal} and Chapters 2 and 3 in \cite{deift1999orthogonal} are introductory references containing these results.

In order to establish context and notation, in Section 2.2 we describe the Lanczos algorithm and its interpretation in terms of orthogonal polynomials.  Some standard references for this matter are Chapter 6 in \cite{trefethen1997numerical} and 
Chapter 6 in \cite{saad2011numerical}. 

In Section 2.3 we introduce the framework in which this paper is developed and formally state the main contributions of our work.

In this paper we use the following notation.  
We use $\rightarrow_P$ to denote convergence in probability.  
For a sequence of events $E_n$, we say $E_n$ occurs with \emph{overwhelming probability} if $\textbf{P}[E_n] \ge 1 - C\exp\{-n^c\}$ for some $c, C > 0$.  For an $n \times n$ matrix $A$ with eigenvalues $\lambda_i$, we say that the \emph{empirical spectral distribution} of $A$ is the atomic probability measure $\frac{1}{n}\sum_{i=1}^n \delta_{\lambda_i}$, where $\delta_x$ denotes the Dirac mass at $x$.  We also let $\Vert A \Vert_{{}}$ denote the spectral norm of $A$ (i.e. the $\ell^2 \to \ell^2$ operator norm).  The notation $\Vert \cdot \Vert$ applied to a vector will refer to the standard Euclidean norm.  We will let $\mathbb{S}^{n-1}$ denote the unit sphere in $\mathbb{R}^n$ and denote the uniform probability measure on $\mathbb{S}^{n-1}$ by $\mathrm{Unif}(\mathbb{S}^{n-1})$. Finally, we will let $\mathrm{Kol}(\cdot,\cdot)$ denote the Kolmogorov distance between two measures, namely,
\[\mathrm{Kol}(\mu, \nu) := \sup_{t \in \mathbb{R}} \left| \mu((-\infty, t]) - \nu((-\infty, t]) \right|. \]
\subsection{Orthogonal polynomials}

For now,  let $\mu$ be a finite Borel measure on $\mathbb{R}$ and assume that its support, which we denote as $\mathrm{supp}(\mu)$, is compact and has infinitely many points. The set of square integrable functions $L^2(\mathbb{R}, d\mu)$ becomes a Hilbert space when endowed with the inner product
 $$\langle f, g\rangle = \int_\mathbb{R} f (x) g(x) d \mu (x) .$$
The hypothesis that $|\mathrm{supp}(\mu)| = \infty$ implies that the monomials $\{1, x, x^2, \dots\}$ are linearly independent in $L^2(\mathbb{R}, d\mu)$.  Hence, we can use the Gram-Schmidt procedure to obtain an infinite sequence of polynomials $p_k(x)$ with $\mathrm{deg}(p_k(x) )= k$ and 
$$\int p_k(x) p_l(x) d \mu(x) = \delta_{kl}. $$

The leading coefficient of $p_k(x)$ is a quantity of interest in this paper and will be denoted by $\gamma_k$. We will denote the monic orthogonal polynomials by $\pi_k(x)$. That is, $\pi_k(x)= \gamma_k^{-1} p_k(x)$ and clearly 

\begin{equation} \label{leadingcoeff} \gamma_k = \left( \int_\mathbb{R} \pi_k^2(x) d \mu(x)  \right)^{-\frac{1}{2}}. \end{equation}

Since $\pi_k(x)$ is orthogonal to all polynomials with degree less than $k$, the polynomial $x^k-\pi_k(x)$ is the orthogonal projection of $x^k$ onto the span of $\{1, \dots, x^{k-1}\}$.  Hence,

$$\int_\mathbb{R} \pi_k^2(x) d\mu(x) = \min_{q\in \Gamma_k} \int_\mathbb{R} q^2(x) d\mu(x),$$
where $\Gamma_k$ denotes the space of monic polynomials of degree $k$. 

Favard's theorem ensures that there is a sequence of real numbers $\alpha_k$ and a sequence of positive real numbers $\beta_k$ such that the following \emph{three-term recurrence} holds:
\begin{align*}
xp_k(x) = \beta_{k-1} p_{k-1}(x)+ \alpha_k p_k(x) +\beta_{k} p_{k+1}(x), \hspace{.3cm} k\geq 1, \\
\text{and} \hspace{.3cm} xp_0(x) = \alpha_0 p_0(x) + \beta_0 p_1(x), \hspace{.3cm} k =0. 
\end{align*}
It is clear from the three-term recurrence that the following identity holds:
\begin{equation}
\label{equationforgamma}
\gamma_k = \left( \prod_{i=0}^{k-1} \beta_i \right)^{-1}.
\end{equation}
These so-called  \emph{Jacobi coefficients} $\alpha_k$ and $\beta_k$ encode all the information of the measure $\mu$. In fact, since the Stieltjes transform of $\mu$ has a continued fraction expansion in terms of its Jacobi coefficients, knowing the few first elements in these sequences allows one to approximate the measure. See Chapter 4.3 in \cite{deift1999orthogonal} for an example. 

We denote by $J_k$ the $k\times k$ Jacobi matrix of $\mu$; that is, $J_k$ is the tridiagonal symmetric matrix with $(J_k)_{ii} = \alpha_{i-1}$ and $(J_k)_{i+1, i} = (J_k)_{i, i+1} = \beta_{i-1}$.  It is a standard fact that $\pi_k(x) = \det(xI-J_k)$ and that in particular, the roots of $p_k(x)$ are exactly the eigenvalues of $J_k$, which are real since $J_k$ is symmetric.  

Another object of importance in this theory is the Hankel matrix of a measure. We will denote $M_k$ the $(k+1)\times (k+1)$ Hankel matrix of $\mu$; in other words, if $m_i$ denotes the $i$th moment of $\mu$, then $(M_k)_{ij} = m_{i+j-2}$ for every $1\leq i, j\leq k+1$. From the elementary theory it is known (see \cite{deift1999orthogonal}, Section 3.1) that if we define $D_k=\det M_k$, then 
\begin{equation}
\label{eqorthpolidentities}
\beta_k = \frac{\sqrt{D_{k-1}D_{k+1}}}{D_k} \hspace{.3cm} \text{and} \hspace{.3cm} \gamma_k = \sqrt{\frac{D_{k-1}}{D_k}}, \hspace{.3cm} k \geq 0,
\end{equation}
where we define $D_{-1} = 1$. Note that the second identity in (\ref{eqorthpolidentities}) implies 
\begin{equation}
\label{eqothgamma}
D_k = \prod_{i=0}^k \gamma_i^{-2}. 
\end{equation}
Moreover, if $\tilde{M}_k(x)$ denotes the matrix obtained by replacing the last row of $M_k$ by the row $(1\: x \: x^2 \: \cdots\:  x^{k})$,  we have the following useful identity:
\begin{equation}
\label{eqdeterminantalid}
p_k(x) = \frac{\det\tilde{M}_k(x) }{\sqrt{D_{k-1} D_k}} .
\end{equation}
Note that in the case in which $\mathrm{supp}(\mu)$ has $n$ points, for $n$ a positive integer, the set of monomials $\{1, x, x^2, \dots \}$ is not linearly independent in $L^2(\mathbb{R}, d\mu)$. Moreover, the Gram-Schmidt procedure stops after $n$ iterations, and hence it only makes sense to talk about the orthogonal polynomials $p_k(x)$ for $k\leq n-1$. However, sometimes it is convenient to define the $n$th monic orthogonal polynomial as the unique monic polynomial of degree $n$ whose roots are the elements of $\mathrm{supp}(\mu)$.  In this case, the facts mentioned previously still hold for $k\leq n$.

\subsection{The Lanczos algorithm}

We understand the Lanczos algorithm as a randomized procedure that takes three inputs: an $n\times n$ Hermitian matrix $A$, a random vector $u$ distributed uniformly in $\mathbb{S}^{n-1}$, and an integer $1\leq k \leq n$. Then, the procedure  outputs a $k\times k$ symmetric  tridiagonal matrix $J_k$ whose diagonal entries will be denoted by $\alpha_i$ for $i=0, \dots, k-1$ and whose subdiagonal and superdiagonal entries will be denoted by $\beta_i$, for $i=0, \dots, k-2$. The eigenvalues of $J_k$ are called the Ritz values and we will usually denote them as $r_1\geq \cdots \geq r_k$. On the other hand, the eigenvectors of $J_k$ give rise (after an orthonormal change of basis determined by the $v_j$ below) to the Ritz vectors, that is, the approximations for the eigenvectors of $A$.   Algorithm \ref{alg:lanczos} below describes how the procedure generates the Jacobi coefficients $\alpha_i$ and $\beta_i$.


        

\begin{algorithm}
\caption{The Lanczos algorithm}
\label{alg:lanczos}

\begin{algorithmic}
\STATE{\textbf{input:} $A$, $k$, $u$}
\STATE{\textbf{initialize: }$v_0=u$}
\FOR{$j = 0, \dots, k-1$}
\STATE{$W_j = \mathrm{span} \{v_0, \dots, v_{j} \}$}
\STATE{$\alpha_{j} = \langle Av_j, v_j \rangle $}
\STATE{$\beta_j = \norm{\mathrm{Proj}_{W_j^{\perp}}(Av_j)}_2$}
\IF{$\beta_j=0$}
\STATE{\textbf{stop}}
\ELSE
\STATE{$v_{j+1} = \frac{\mathrm{Proj}_{W_j^{\perp}}(Av_j)}{\norm{  \mathrm{Proj}_{W_j^{\perp}}(Av_j)}_2}$}
\ENDIF
\ENDFOR
\RETURN $J_k$
\end{algorithmic}
\end{algorithm}

This algorithm has a natural interpretation in terms of orthogonal polynomials.  To every $u \in \mathbb{S}^{n-1}$ we can associate a measure supported on the spectrum of $A$ as follows. Let $\lambda_1 \geq \cdots \geq \lambda_n$ be the eigenvalues of $A$ and $u_1, \dots, u_n$ be the coordinates of $u$ when writen in the eigenbasis of $A$. We define the probability measure
\begin{equation}
\label{eqmuu}
\mu^u = \sum_{i=1}^n u_i^2 \delta_{\lambda_i}. 
\end{equation}

In the language of functional analysis, $\mu^u$ is the spectral measure of the operator $A$ induced by the vector state $u$; that is, $\langle f(A)u, u \rangle  = \int f(x)\,d\mu^u(x)$ for all (say) polynomials $f$. Note that the expectation of the random measure $\mu^u$ is just the empirical spectral distribution of $A$, namely,
$$\frac{1}{n}\sum_{i=1}^n \delta_{\lambda_i}.$$

It is not hard to see that if $p_j(x)$ are the orthogonal polynomials with respect to $\mu^u$, then 
$v_j = p_j(A) u$. Hence, the coefficients $\alpha_j$ and $\beta_j$ outputed by the Lanczos algorithm are the Jacobi coefficients of the measure $\mu^u$, and the Ritz values after $k$ iterations are the roots of $p_k(x)$. 

As a last remark, observe that the output of Algorithm \ref{alg:lanczos} scales linearly with $A$. Hence, to simplify notation, in some of the proofs below we will start by assuming that $\norm{A}_{{}} =1$.

\subsection{Statement of results}  In the remainder of the paper the input matrix, which is assumed to be Hermitian, will be fixed and denoted by $A$.  We will use $n$ to denote the dimension of $A$ and usually the number of iterations of the procedure will be denoted by $k$.  Note that the Jacobi coefficients $\alpha_i$ and $\beta_i$ are assigned during the $i$th iteration of the algorithm and are unchanged during future iterations.

Since for our analysis it is necessary to compare outputs of the algorithm resulting from different input vectors $u\in \mathbb{S}^{n-1}$, we will stress this dependence by viewing the respective quantities as a function of $u$ and denoting them by $\alpha_i(u), \beta_i(u), r_i(u), \gamma_i (u)$, $p_k^u(x)$, $v_i(u)$, and $J_k(u)$. Depending on the context, these quantities will also be thought of as random variables, random polynomials, random vectors, and random matrices, respectively. 

Most of our results make a technical assumption, which we have named \emph{equidistribution}, about the geometry of the spectrum of $A$. We define and motivate this concept below. 
\bigskip

\noindent \textbf{Equidistribution.} The behaviour of the Lanczos algorithm depends on the spectrum of the input matrix $A$, and what might be true for typical matrices can fail for particular choices of $A$. It is hence challenging, when trying to obtain general theoretical statements,  to identify a quantifiable feature of the spectrum that can serve as an assumption but is not too restrictive. This has been done in different ways in previous work, most notably in the seminal work of Saad \cite{saad1980rates}---see \cite[Section 6.6]{saad2011numerical} for a succinct exposition.  For example, many of the results in \cite{saad1980rates} are stated in terms of the  quantities
$$t_i^{(k)} = \min_{p \in \mathcal{P}_{k-1}^{(i)}} \max_{j: j\neq i} |p(\lambda_j)|,$$
where $\lambda_1\geq \cdots \geq \lambda_n $ are the eigenvalues of $A$ and $\mathcal{P}_{k-1}^{(i)}$ denotes the set of all polynomials of degree not exceeding $k-1$ and satisfying $p(\lambda_i) =1$. In the present work we use a more geometric notion to state our results.

\begin{definition}[Equidistribution]\label{defnequidistribution}
 Let $\Lambda$ be any finite set of $n$ real numbers.  Let $\delta$ and $\omega$ be positive real numbers and let $j$ be a natural number.  We say that $\Lambda$ is $(\delta, \omega, j)$-\emph{equidistributed} if for any finite set $T$ of at most $j$ real numbers it holds that

\[ \left| \left\{ \lambda \in \Lambda : \frac{1}{|T|} \sum_{t \in T} \log |\lambda - t| \ge \log \omega \right\} \right|  \geq \delta n.\]
\end{definition}

  Intuitively, the spectrum is equidistributed if it is not grouped in a small number of tight clusters (see Examples \ref{exampleparameters} and \ref{exampleparamclusters} below). As we will show in Section 4, the family of well equidistributed point sets includes, but is not limited to, those sets obtained by discretizing an absolutely continuous distribution. 


\begin{example}
\label{exampleparameters}
Let $\Lambda$ be the set of $n$ equally spaced points from $1/n$ to $1$, inclusive.  This represents a discretization of the uniform measure $\mu = \mathrm{Unif}([0,1])$. In Section 4.1 we will show that for $j\leq \frac{n}{16}$, the set $\Lambda$ is $(\delta, \omega, j)$-equidistributed for $\delta = 1/4$ and $\omega = 4e^{-2}$. 
\end{example}

\begin{example}
\label{exampleparamclusters}
Now consider a set (or multiset) $\Lambda$   of $n>0$ points grouped in $m$ equally spaced small clusters. To make this precise, fix  two parameters $ \varepsilon, g >0$ and consider $-1= a_1 \leq b_1 < a_2\leq b_2 < \cdots < a_m \leq b_m =1$ such that for every $i= 1, \dots, m$ we have $b_i-a_i = \varepsilon$ and $a_{i+1}-b_i = g$. We think of $\varepsilon$ as  \emph{small} with respect to $g$ and of $m$ as \emph{small} with respect to $n$.  If $\Lambda \subset \bigcup_{i=1}^m [a_i, b_i]$ with $|\Lambda \cap [a_i, b_i] | \geq \lfloor \frac{n}{m} \rfloor$ for every $i=1, \dots, m$, then $\Lambda$ is $(\frac{m-j}{m}, g, j)$-equidistributed and $g \approx 2/m$. 

Note that in this case we have good equidistribution parameters unless $j\approx  m$. In Section 4  we  give a generalization of this assertion in Observation \ref{obsboringcase}. 
\end{example}

\bigskip

\noindent \textbf{Concentration of the output.} We state our main  result about concentration of the first few Jacobi coefficients in terms of the equidistribution parameters $\delta, \omega$. 

\begin{theorem}[Concentration of Jacobi coefficients after $i$ iterations]\label{thmhessconc}
Suppose the spectrum of $A$ is $(\delta, \omega, i)$-equidistributed for some $\delta, \omega > 0$ and $i \in \mathbb{N}$.  Let $\tilde{\alpha}_i$ and $\tilde{\beta}_i$ denote the medians of the Jacobi coefficients $\alpha_i(u)$ and $\beta_i(u)$, respectively. Then for all $t > 0$, the probabilities
$\mathbf{P}[|\alpha_{i}(u) - \Tilde{\alpha_{i}}| > t\Vert A \Vert]$ and $\mathbf{P}[|\beta_{i}(u) - \Tilde{\beta_{i}}| > t\Vert A \Vert]$ are both bounded above by 
\begin{equation}
\label{ineqConcentration}
2\exp\left\{-\frac{\min\{\delta, 1/50\}^2}{32}n\right\} + 2 \exp \left\{- \frac{1}{64}\left(\frac{\omega}{4 \Vert A \Vert} \right)^{2i}\delta^2 t^2 n\right\}. 
\end{equation}

\end{theorem}

To clarify the advantages and limitations of the above result we include some remarks. 

\begin{remark} 
The equidistribution parameters $\delta, \omega$ appearing in the above theorem are typically quite moderate in magnitude and are easy to compute if one can obtain explicit bounds for certain integrals with respect to the spectral distribution of $A$. Note that $\omega \le \Vert A\Vert$ (by taking $T = \{0\}$ in Definition \ref{defnequidistribution}) and that $\omega$ scales linearly with $A$.  As a result, $\omega/\Vert A \Vert$ is typically of constant size independent of $n$ in applications. Since $\omega/\Vert A \Vert < 1$, Theorem \ref{thmhessconc} yields concentration for $i$ at most logarithmic in $n$. Besides Examples \ref{exampleparameters} and \ref{exampleparamclusters}  above, in Section 4.1 we give more examples and a  detailed discussion on how to compute these parameters. 
\end{remark}

\begin{remark}
Ultimately because the Jacobi coefficients $\alpha_i(u), \beta_i(u)$ turn out to be Lipschitz only on a proper subset of the sphere, we only obtain concentration of the coefficients about their medians, not their means.
\end{remark}

\begin{remark}\label{pseudospectrum}
Using the same techniques one can prove a result analogous to Theorem \ref{thmhessconc} in the case where $A$ is not Hermitian, and even not normal.  In the non-Hermitian case, the Lanczos algorithm is called the \emph{Arnoldi algorithm} and has similar applications to those of the Lanczos algorithm.  
\end{remark}


In Section 4 we show how Theorem \ref{thmhessconc} can be used to obtain concentration, in the same regime of number of iterations, of the Ritz values and Ritz vectors.   See Propositions \ref{propConcRitzVals} and  \ref{propConcRitzeigenvectors} for  precise statements.


The rest of the paper focuses on studying the location of the outputs of the Lanczos algorithm. 
\bigskip

\noindent \textbf{Undetected outliers of the spectrum.} In Section 5.1, we show that if $k$ is a certain fraction of $\log n$, the Ritz values obtained after $k$ iterations are contained in a small blow-up of the convex hull of the bulk of the spectrum of $A$. This complements classical guarantees which show that for some multiple of $\log n$, say, $K$, the Lanczos algorithm approximates with high accuracy the outliers of the spectrum of $A$ when $K$ iterations are performed.  Our results are quantitative and use our notion of equidistribution. 

\begin{theorem}\label{thm:nonasylocation}
Suppose the spectrum of $A$ is $(\delta, \omega, j)$-equidistributed  for some $\delta, \omega > 0$ and $j \in \mathbb{N}$.  Let $M$ be the diameter of the spectrum of $A$. Let $R$ be a real number and let $0 < c < 1/2$, and suppose there are at most $m\le \min\{0.02n, 2 n^\alpha\}$ ``outliers,'' eigenvalues of $A$ lying above $R$, for some $\alpha < 1-c$.  Let $g = \max_{1 \le i \le n}\{\lambda_i - R\} $ and let $\kappa > 0$.  Then for up to \[k = \min\left\{j, \frac{1}{2 \log \frac{M}{\omega}}\left(c \log n + \log \frac{\kappa \delta}{2 m g}\right)  \right\}\]
iterations, the probability that the top Ritz value exceeds $R + \kappa$ is at most 
\[2\exp\left\{-\frac{\min\{\delta,  1/50\}^2}{32} n\right\} + 2\exp \left\{ -\frac{1}{16} n^{1-2c} \right\}\] 
for $n > e^{\frac{1}{1-c-\alpha}}$.  
\end{theorem}
The strength of the above result might be obscured by the appearance of several unintuitive parameters.  For the reader's benefit we include an example below, and to provide a slightly different perspective, we include an asymptotic version of the above result, namely, Proposition \ref{prop:outliers}. 

\begin{example}
\label{exampleparameters2}
Let $n>0$ and let $A$ be a matrix whose spectrum consists of $n-1$ equally spaced points from $2/n$ to $1$ inclusive, together with an outlier of value 1.1 (compare with Figure \ref{fig:outlier}). In Section 4.1 we will show that for $j \leq n/16$ the spectrum of $A$ is $(1/4, 4e^{-2}, j)$-equidistributed.  

In order to apply Theorem \ref{thm:nonasylocation}, we also note that in this case $M = 1.08$, $m=1$, and $g=10^{-1}$.  Take $\kappa = 10^{-4}$. Then, for any $0 < c< 1/2$, the Ritz values of the Lanczos algorithm on $A$ after $\lfloor \frac{7c}{10}\log n -7/2\rfloor$ iterations will  be contained in the interval $[2/n, 1+10^{-4}]$ with  overwhelming probability. 
\end{example}


\begin{proposition} \label{prop:outliers}
Let $(A_n)_{n=1}^\infty$ be a sequence of $n \times n$ Hermitian matrices with uniformly bounded norm. Assume their empirical spectral distributions $\mu_n$ converge in distribution to a measure $\mu$ with nontrivial absolutely continuous part, and further assume $\text{Kol}(\mu_n, \mu) = O(1/\log n)$.  Suppose there exists $m \in \mathbb{N}$ such that each $A_n$ has at most $m$ eigenvalues (``outliers'') greater than $R$, where $R$ denotes the right edge of the support of $\mu$.  

Then there exists $c > 0$ such that for every $\kappa > 0$, the Ritz values of Lanczos applied to $A_n$ after $c \log n$ iterations are bounded above by $R+\kappa$ with overwhelming probability for $n$ sufficiently large (depending on how small the gap $\kappa$ is chosen.)
\end{proposition}

\begin{remark}
Suppose all eigenvalues of $A$ lie in $[-1, 1]$ except for one outlier $\lambda_1 = 1+\ep$.  Suppose we wish for a Ritz vector $p(A)u = \sum_{i=1}^n p(\lambda_i) u_i$ to approximate $u_1$, the top eigenvector.  Then $p(\lambda_1)$ must be polynomially larger than $p(\lambda_2), \dots, p(\lambda_n)$, since $|u_i| \sim 1/\sqrt{n}$ for all $1 \le i \le n$ with high probability.  If one further imposes that $p(\lambda_1)$ exceeds $p$ by a polynomial factor on the entire interval $[-1, 1]$, then the degree of $p$ (and hence the number of iterations $k$) must be at least $O(\log n / \sqrt{\ep}),$ by the Markov brothers' inequality and properties of Chebyshev polynomials.  However, if the $\lambda_i$ are in $j$ tight clusters for $j$ small, one could certainly have $p(\lambda_2), \dots, p(\lambda_n)$ small and $p(\lambda_1)$ large for some $p$ of degree $j$.  Thus, an assumption like equidistribution is natural for Theorem \ref{thm:nonasylocation}.
\end{remark}

\bigskip

\noindent \textbf{Location of Jacobi coefficients and Ritz values } We will now  work on the setting of Proposition \ref{prop:outliers}, that is, we consider a probability measure $\mu$ and a sequence of matrices $A_n$ whose empirical spectral distributions $\mu_n$ converge to $\mu$. We give a result about the locations of the Ritz values and Jacobi coefficients when at most $d \sqrt{\log n}$ iterations are performed, with $d$ depending only on $\mu$ and the speed of convergence of the sequence $\mu_n$. Essentially, we show that in this regime the Jacobi matrix after $k$ iterations is sharply concentrated around the $k$th Jacobi matrix of the measure $\mu$. 

\begin{theorem}[Location of Jacobi coefficients]
\label{Thmsqrtiter} Let $(A_n)_{n=1}^\infty$ be a sequence of $n \times n$ Hermitian matrices with uniformly bounded operator norm. Assume their empirical spectral distributions $\mu_n$ converge in distribution to a measure $\mu$ with nontrivial absolutely continuous part, and further assume $\text{Kol}(\mu_n, \mu) = O(n^{-c})$ for some $c > 0$.

 Then there is a constant $d>0$ dependent on $\mu$ and $c$, such that  for any sequence of integers $1\leq k_n\leq d\sqrt{\log n} $ we have 
$$\norm{J_{k_n}(u)- J_{k_n}(\mu)}_{{}} \longrightarrow_P 0, $$
where $J_{k_n}(u)$ denotes the Jacobi matrix output by the Lanczos algorithm applied to $A_n$ under the input $u\sim \mathrm{Unif}(\mathbb{S}^{n-1})$ after $k_n$ iterations,  and where $J_{k_n}(\mu)$ is the $k_n$-th Jacobi matrix of the measure $\mu$.
\end{theorem}

Note that Theorem \ref{Thmsqrtiter} may be of particular relevance in applications where an infinite-dimensional operator is discretized with the goal of computing its density. In essence, Theorem \ref{Thmsqrtiter} states that, in this situation, the first iterations of the Lanczos algorithm are an accurate approximation of the true Jacobi coefficients of the measure $\mu$, and hence the procedure gives valuable information to recover the limiting  measure. 

From the above proposition, a standard application of the Weyl eigenvalue perturbation inequality  yields the following proposition. 

\begin{proposition}[Location of the Ritz values]
\label{cormain}
Using the same notation as in Theorem \ref{Thmsqrtiter}, let $\vec{r}_{k_n}(u) = (r_1(u), \dots, r_{k_n}(u))$, where $r_1(u)\geq \cdots \geq r_{k_n}(u)$ are the random Ritz values of the Lanczos algorithm after $k_n$ iterations are performed. Then under the assumptions in Theorem \ref{Thmsqrtiter}, we have that 
$$\norm{\vec{r}_{k_n}(u)- \vec{r}_{k_n}(\mu)}_{L^\infty(\mathbb{R}^{k_n})} \longrightarrow_P 0,$$
where $\vec{r}_{k_n}(\mu)$ is the vector whose entries are the roots of the $k_n$-th orthogonal polynomial with respect to $\mu$ in decreasing order.  
\end{proposition}

It remains an open question if similar results can be obtained when $O(\log n)$ iterations are performed. See Section 6 for open questions and  further research. 

\begin{remark} \label{rem:fixedk}
For \emph{fixed} $k$, \cite[Theorem 4]{gautschi1968construction} states that in the deterministic setting of a weakly convergent sequence of measures $\mu_n \to \mu$ with each $\mu_n$ supported on $n$ points, the Jacobi coefficients $\alpha_i, \beta_i$ of $\mu_n$ for $i \le k$  converge to those for $\mu$ (and therefore, the same holds for the $i$th orthogonal polynomial and its roots.)  To apply this to the Lanczos algorithm, one must use the fact that for independent initial vectors $u_n \in \mathbb{S}^{n-1}$, almost surely the measures $\mu_n^{u_n}$ (defined in \ref{eqmuu}) converge weakly to $\mu$.  One therefore obtains the same convergence as in the deterministic setting, almost surely, for the random $\alpha_i, \beta_i$ output by Lanczos.  To make the convergence quantitative, one must assume something about the rate of convergence of the $\mu_n$ in the hypothesis; this is what is done in Theorem \ref{Thmsqrtiter}.  Our result holds in the wider range $k = O(\sqrt{\log n})$.  We also trade almost-sure convergence for a polynomial rate of convergence holding with overwhelming probability; this rate appears in the proof.
\end{remark}

\section{Applying the local L\'evy lemma}

\subsection{Strategy}

The well known L\'evy lemma states, in a quantitative way, that if $f: \mathbb{S}^{n-1}\to \mathbb{R}$ is a Lipschitz function, then $f(u)$ is a random variable concentrated around its median. See  Chapter 5.1 in \cite{vershynin2018high} for a detailed discussion. In this direction,  the main obstacle for showing concentration of the random variables $\alpha_i(u)$ and $\beta_i(u)$ is that the functions $\alpha_i, \beta_i: \mathbb{S}^{n-1}\to \mathbb{R}$ are not Lipschitz on the entire sphere. However, we will be able to show that these functions are  Lipschitz in a large region of the sphere, which is a common idea in geometric functional analysis. We will use a local version of L\'evy's lemma, which is recorded as Corollary 5.35 in \cite{aubrun2017alice}, and which we restate below with explicit universal constants. 

\begin{lemma}[Local L\'evy lemma]\label{lemmalevy}
 Let $\Omega \subset \mathbb{S}^{n-1}$ be a subset of measure larger than $3/4$. Let $f : \mathbb{S}^{n-1} \to \mathbb{R}$ be a function such that the restriction of f to $\Omega$ is Lipschitz with constant $L$ (with respect to the geodesic metric on the sphere). Then, for every $\varepsilon > 0$,
$$\mathbf{P}[|f(u) - \tilde{f} | > \varepsilon] \leq  \mathbf{P}[u \in \mathbb{S}^{n-1}\setminus \Omega] + 2 \exp\{-4n \varepsilon^2/L^2\},$$
where $\tilde{f}$ is the median of $f(u)$ and where $u\sim \mathbb{S}^{n-1}$. 
\end{lemma}

As the function $f$ is allowed to blow up on the small subset $\mathbb{S}^{n-1} \setminus \Omega$, one cannot expect a similar result to hold for concentration around the mean.  

In order to identify the correct region of the sphere in which the functions $\alpha_{i}$ and $\beta_i$ are Lipschitz, we need a local version of the notion of Lipschitz constant. In what might be a slight departure from standard definitions, we will define \emph{local Lipschitz continuity} as follows.

\begin{definition}
\label{deflocalLipschitzness}
Let $(X_1, d_1)$ and $(X_2, d_2)$ be metric spaces. A function $f: X_1 \to X_2$ is said to be locally Lipschitz continuous  with constant $c$ at $x_0 \in X_1$ if for every $c' > c$ there is a neighborhood $U\subset X_1$ of $x_0$ such that 
$$d_2(f(x),f(y)) \leq c' d_1(x, y) \hspace{.3cm} \forall x , y \in U.$$
\end{definition}

\begin{remark}
For $f$ defined on an open subset of $\mathbb{S}^{n-1}$, we have that $f$ is locally Lipschitz continuous with constant $c$ with respect to the geodesic metric if and only if it is locally Lipschitz continuous with the same constant with respect to the Euclidean (``chordal'') metric.
\end{remark}
It is obvious that if a function is locally Lipschitz with constant $c$ on every point of a convex set, then the function is globally Lipschitz on the set with the same constant $c$. However, if the convexity assumption is dropped, a similar conclusion is not guaranteed in general and in order to obtain a global Lipschitz constant the geometry of the set should be analyzed.

\begin{definition}
Let $K > 0$ and  $(X, d)$ be a metric space. We say that $S_1\subset X$  is $K$-connected in $S_2$  with $S_1 \subset S_2 \subset X$ if for every $x, y \in S_1$ there is a rectifiable Jordan arc $\alpha: [0, 1] \to S_2$ with  $\alpha(0) = x$ and $\alpha(1) = y$, such that the length of the trace of  $\alpha$ is less than or equal to $K d(x, y)$.  
\end{definition} 

Now that we have introduced the notion of  $K$-connected set we can generalize what we observed for convex sets.  

\begin{lemma}
\label{lemmalocalglobal}
Let $(X_1, d_1)$ and $(X_2, d_2)$ be metric spaces. Assume that $S_1 \subset X_1$ is $K$-connected in $S_2 \subset X_1$ and let   $f: X_1 \to X_2$ satisfy that for every $x_0 \in S_2$,  $f$ is locally Lipschitz at $x_0$ with constant $c$. Then $f$ is globally Lipschitz on $S_1$ with constant $cK$.    
\end{lemma} 

\begin{proof}
Fix $x, y \in S_1$ and $\varepsilon > 0$. We will show that $d_2(f(x),f(y)) \leq (c+ \varepsilon) K d_1(x, y) $. Consider a rectifiable Jordan arc $\alpha: [0,1] \to X_1$, such that $\alpha(0)=x$, $\alpha(1)=y$, $\alpha([0,1])\subset S_2$ and the length of $\alpha$ is at most $K d_1(x, y)$. 

Since the trace of $\alpha$ is contained in $S_2$, for every $w\in \alpha([0,1])$ we can take an open ball $U_w$ containing $w$ such that $f$ is $(c+\varepsilon)$-Lipschitz on $U_w$. Moreover, observe that since $\alpha$ is continuous and injective, for every $w \in \alpha([0,1])$ we can take $U_w$ small enough such that $\alpha^{-1}(U_w)$ is connected and hence an open interval in $[0,1]$.

By compactness  of $\alpha([0,1])$  we may take $w_1, \dots, w_n \in  \alpha([0,1])$ such that $\{U_{w_i}\}_{i=1}^n$ is a minimal cover for $\alpha([0,1])$. Now, since each $\alpha^{-1}(U_{w_i})$ is connected, and the cover is minimal, we have that $\alpha^{-1}(U_{w_i}) \cap \alpha^{-1}(U_{w_{i+1}}) \neq \emptyset$ for every $ 1\dots, n-1$. 

Furthermore, we will now see that we can modify the sequence of $w_i$ such that $w_{i+1}\in U_{w_i}$ for every $i=1, \dots, n-1$. Assume that this does not hold and let $i$ be the smallest index for which $w_{i+1}\notin U_{w_{i}}$.  Now take some $t \in \alpha^{-1}(U_{w_i}) \cap \alpha^{-1}(U_{w_{i+1}})$ and define $w' = \alpha(t)$. We construct a new sequence $\tilde{w}_1, \dots, \tilde{w}_{n+1}\in \alpha([0,1])$ by taking $\tilde{w}_j = w_{j}$ for $j <i$, $\tilde{w}_i = w'$,  $\tilde{w}_{j+1} = w_j$ for $j\geq i$, and  $U_{\tilde{w}_i}$ to be equal to $U_{w_{i+1}}$. Observe that for the new sequence of points $(\tilde{w}_i)_{i=1}^{n+1}$ in $\alpha([0,1])$ and sequence  of open balls $U_{\tilde{w}_i}$ it holds that $\tilde{w}_{j+1} \in U_{\tilde{w}_j}$ for all $j\leq i$. By iterating this process we will obtain a finite sequence with the desired property. So, in what follows we can assume without loss of generality that $w_{i+1}\in U_{w_i}$ for every  $i=1, \dots, n-1$.  We then will have 
$$d_2(f(w_i), f(w_{i+1})) \leq (c+ \varepsilon) d_1(w_i, w_i+1).$$
Using the triangle inequality and the fact that $\sum_i d_1(w_i, w_{i+1})$ is bounded by the length of the trace of $\alpha$ the result follows. 
\end{proof}

In the following section the local Lipschitz constants of the functions $\alpha_{i}(u)$ and $\beta_i(u)$ are shown to be related to the orthogonal polynomials of the measure $\mu^u$.

\subsection{Local Lipschitz constants for Jacobi coefficients}
  As can be seen from Algorithm \ref{alg:lanczos}, the dependence of the quantities $\alpha_{i}(u)$, $\beta_i(u)$, and $v_j(u)$ on $u$ is highly nonlinear, which makes it complicated to show that such quantities are stable under perturbations of the input vector $u$. Here we exploit the fact that during  every iteration of the Lanczos algorithm only locally Lipschitz operations are performed. The analysis of the compound effect of iterating the procedure yields a bound on the local Lipschitz constant of the quantities of interests. This bound is exponential in the number of iterations, which is enough to obtain concentration results when  $O(\log (n))$ iterations are performed. In what follows, recall that $\gamma_i(u)$ denotes the leading coefficient of the $i$th orthonormal polynomial with respect to the measure $\mu^u$ defined in (\ref{eqmuu}).

\begin{proposition}
\label{propositionLipsconstant}
Fix $\tilde{u} \in \mathbb{S}^{n-1}$ and let $v_j(u)$ be as in Algorithm \ref{alg:lanczos}.   Then, for any $0\leq j \leq n-1$, the functions $v_j(u)$ are locally Lipschitz at $\tilde{u}$ with constant $(4 \norm{A}_{{}})^j \gamma_{j}(\tilde{u})$.
\end{proposition}

\begin{proof}

We proceed by induction. For $j=0$, recall $v_0(u) = u$ and  $\gamma_{0}(\tilde{u}) =1$; the statement follows. Now assume the proposition is true for some $j\geq 0$.
 For every $x\in \mathbb{S}^{n-1}$  denote $W_{x} = \mathrm{span}\{ v_0(x) = x, v_1(x), \dots, v_j(x)\}$ and  for any subspace $W \leq \mathbb{R}^n$ by $\mathrm{Proj}_W$ we mean the orthogonal projection onto  $W$.

  Take $x, y \in \mathbb{S}^{n-1}$ in a neighborhood $\mathcal{U}$ of $\tilde{u}$ to be determined and note that
\begin{align}
\label{ineqproj1}
& \norm{\mathrm{Proj}_{W_x^\perp}(Av_j(x))-\mathrm{Proj}_{W_y^{\perp}}(A v_j(y))} \nonumber \\ & \leq \norm{\mathrm{Proj}_{W_x^\perp}(A (v_j(x)-v_j(y)))} + \norm{(\mathrm{Proj}_{W_x^\perp}-\mathrm{Proj}_{Wy^\perp})(Av_j(y))}  \nonumber 
\\ & = \norm{\mathrm{Proj}_{W_x^\perp}(A (v_j(x)-v_j(y)))} + \norm{(\mathrm{Proj}_{W_x}-\mathrm{Proj}_{W_y})(Av_j(y))}. 
\end{align}

From the induction hypothesis we have that, for any $\varepsilon> 0$, we can choose $\mathcal{U}$ small enough so that 
\begin{equation}
\label{ineqproj2}
\norm{\mathrm{Proj}_{W_x^\perp}(A (v_j(x)-v_j(y)))} \leq \norm{A}_{{}} \norm{v_j(x)-v_j(y)} \leq \norm{A}_{{}}((4\norm{A}_{{}})^j \gamma_{j}(\tilde{u})+\varepsilon)\norm{x-y}.     
\end{equation}

On the other hand, from Algorithm \ref{alg:lanczos} it follows that $\beta_{i}(\tilde{u})\leq \norm{A}_{{}}$ for every $i=0, \dots, n-1$, so in view of (\ref{equationforgamma}), the $\norm{A}_{{}}^i \gamma_i(\tilde{u})$ form an increasing sequence.  It then follows that
$$\sum_{i=0}^j (4 \norm{A}_{{}} )^i \gamma_{i}(\tilde{u}) \leq  \sum_{i=0}^j 4^i \norm{A}_{{}}^{j} \gamma_{j}(\tilde{u}) \le \frac{4^{j+1} \norm{A}_{{}}^j \gamma_{j}(\tilde{u})}{3} . $$
 For any unit vector $w$, by the triangle inequality, we have that
\begin{equation}
\label{eqrtg}
 \norm{\mathrm{Proj}_{W_x}(w)- \mathrm{Proj}_{W_y}(w)} \leq \sum_{i=0}^j \norm{\langle v_i(x), w \rangle v_i(x) - \langle v_i(y), w\rangle v_i(y)}
\end{equation}
and we can bound each term on the right-hand side of (\ref{eqrtg}) as follows:
\begin{align*}
\norm{\langle v_i(x), w \rangle v_i(x) - \langle v_i (y), w \rangle v_i(y)} &\leq |\langle v_i(x)-v_i(y), w \rangle| + \norm{v_i(x)-v_i(y) } |\langle v_i(y), w\rangle  | \\ & \leq \norm{v_i(x)-v_i(y)}\norm{w}+ \norm{v_i(x)-v_i(y)} \norm{v_i(y)}\norm{w} \\ & \leq 2 (4\norm{A}_{{}})^i \gamma_{i}(\tilde{u})\norm{x-y}.
\end{align*}
 Hence, adding over $i$ we obtain  
$$\norm{\mathrm{Proj}_{W_x}(w)- \mathrm{Proj}_{W_y}(w)} \leq \frac{2}{3}\cdot 4^{j+1} \norm{A}_{{}}^j \gamma_{j}(\tilde{u}) \norm{x-y},$$
which implies that $\norm{\mathrm{Proj}_{W_x}-\mathrm{Proj}_{W_y}}_{} \leq \frac{2}{3}\cdot 4^{j+1} \norm{A}_{{}}^j \gamma_{j} (\tilde{u}) \norm{x-y}$ and hence 
\begin{equation}
\label{ineqproj3}
\norm{(\mathrm{Proj}_{W_x}-\mathrm{Proj}_{W_y})(A v_j(y))} \leq \frac{2}{3}\cdot (4 \norm{A}_{{}})^{j+1} \gamma_{j} (\tilde{u}) \norm{x-y}.    
\end{equation}
 Putting together inequalities (\ref{ineqproj1}), (\ref{ineqproj2}), and (\ref{ineqproj3}), we get for any $x, y\in \mathcal{U}$  that
$$\norm{\mathrm{Proj}_{W_x^\perp}(Av_j(x))-\mathrm{Proj}_{W_y^\perp}(A v_j(y))}\leq (4\norm{A}_{{}})^{j+1} \gamma_{j}(\tilde{u}) \norm{x-y}. $$
With this we have established that the function $u \mapsto \mathrm{Proj}_{W_u^{\perp}}(Av_j(u))$ is locally Lipschitz at $\tilde{u}$ with constant $ (4\norm{A}_{{}})^{j+1} \gamma_{j}(\tilde{u})$. Now consider the function $f : \mathbb{R}^n \to \mathbb{R}^n$ defined by $f(x) = x/\norm{x}$. It is easy to show that for any $x_0\neq 0$, $f$ is locally Lipschitz at $x_0$ with constant $1/ \norm{x_0}$. Now recall that by definition $\beta_{j}(\tilde{u}) = \norm{ \mathrm{Proj}_{W_{\tilde{u}}^\perp}(Av_j(\tilde{u}))}$. Since the composition of locally Lipschitz functions is locally Lipschitz with the constant being the product of the constants of each of the functions in the composition, we have that the function 
$$u\mapsto v_{j+1}(u) = f(\mathrm{Proj}_{W_u^\perp}(Av_j(u)))$$
is locally Lipschitz at $\tilde{u}$ with constant  $\frac{(4\norm{A}_{{}})^{j+1} \gamma_{j}(\tilde{u})}{\beta_{j}(\tilde{u})} = (4\norm{A}_{{}})^{j+1} \gamma_{j+1}(\tilde{u})$, where this equality follows from  (\ref{equationforgamma}).   
\end{proof}

\begin{proposition}
\label{propositionhessenberglipschitz}
For any $0\leq j  \leq n-1$ and any $\tilde{u}\in \mathbb{S}^{n-1}$, the function $\alpha_{j}(u)$ is locally Lipschitz at $\tilde{u}$ with constant $\frac{1}{2}\cdot (4\norm{A}_{{}})^{j+1} \gamma_{j}(\tilde{u})$, while  $\beta_{j}(u)$ is locally Lipschitz at $\tilde{u}$ with constant $(4\norm{A}_{{}})^{j+1} \gamma_{j}(\tilde{u})$. 
\end{proposition}

\begin{proof}
We will use the same notation as in Proposition \ref{propositionLipsconstant}.  Recall from Algorithm \ref{alg:lanczos} that $\alpha_j(u) = \langle A v_j(u),  v_j(u) \rangle$.  Note that the local Lipschitz constant of the function $u \mapsto A v_j (u)$ is obtained by multiplying the local Lipschitz constant of $v_j(u)$ by $\norm{A}_{{}}$ . Then, for any $\varepsilon$ we can pick  $\mathcal{U}$ to be a small enough neighborhood of $\tilde{u}$ such that for any  $x, y \in \mathcal{U}$ we have  
\begin{align*}
   |\alpha_j(x)- \alpha_j(y)|  &= |\langle Av_j(x),  v_j(x)\rangle- \langle Av_j(y),  v_j(y)\rangle| \\ &\leq |\langle A(v_j(x)-v_j(y)), v_j(x)\rangle| + |\langle  A v_j(y), v_j(x)- v_j(y) \rangle | \\
   &\leq 2\cdot (4^{j} \norm{A}_{{}}^{j+1} \gamma_i (\tilde{u})+\varepsilon)\norm{x-y}.
\end{align*}

 On the other hand, since  $\beta_j(u) =\norm{ \mathrm{Proj}_{W_u^\perp}(Av_j(u)))}$ and we  established in the proof of Proposition \ref{propositionLipsconstant} that this function is locally Lipschitz with constant $(4\norm{A}_{{}})^{j+1}\gamma_{j}(\tilde{u})$, the proof is concluded. 
\end{proof}

\begin{remark}
The local Lipschitz constants presented in the above statements can be improved; the term $4^{j}$ next to $\norm{A}_{{}}^j \gamma_{j}(\tilde{u})$ was chosen for the sake of exposition. Nevertheless, it seems complicated to show that the quantities $v_j(u)$ are locally Lipschitz at $\tilde{u}$ with a constant of the form $C_j \norm{A}_{{}}^{j} \gamma_{j}$ and $C_j$ subexponential. In any case, the term $\norm{A}_{}^j \gamma_{j}$ is typically exponential in $j$, so an improvement on $C_j$ would not yield an asymptotic improvement to the final result if the same level of generality is considered. However, as we point out in Section 6, sharpening our constants is of relevance for applications. 
\end{remark}

\subsection{Incompressibility}
 In Section 4, we will see that our upper bounds for the local Lipschitz constants of the Jacobi coefficients go to infinity if $u$ becomes too close to a sparse vector, roughly speaking.  So we only have a good local Lipschitz constant in a certain region of the unit sphere that avoids sparse vectors.  In order to upgrade our local Lipschitz constant to a global Lipschitz constant, we must prove

\begin{enumerate}[(1)]
    \item that this region is large enough to apply the local L\'evy lemma (Lemma \ref{lemmalevy}) and 
    \item that this region is $K$-connected for a small enough $K$.
\end{enumerate}

First we give this region a name.  Loosely inspired by the compressed sensing literature (see, for example, \cite{vershynin2009role}), we say that a vector $u$ in $\mathbb{S}^{n-1}$ is $(\delta, \ep)$-\emph{incompressible} if each set of at least $\delta n$ coordinates carries at least $\ep$ of its ``$\ell^2$ mass.''  Otherwise, we say that $u$ is $(\delta, \ep)$-\emph{compressible}.  We denote the set of $(\delta, \ep)$-incompressible vectors in $\mathbb{S}^{n-1}$ by $I_n(\delta, \ep)$ and record the formal definition below.

\begin{definition}
\[I_n(\delta, \ep) = \left\{ u \in \mathbb{S}^{n-1} : \sum_{i \in S} u_i^2 > \ep \text{ for all } S \subseteq \{1, 2, \dots, n \}, |S| \ge \delta n \right\}.\]
\end{definition}

For incompressible $u$ we prove an adequate bound on the local Lipschitz constant in Proposition \ref{propgammabound}. Fortunately, a uniform random unit vector $u$ is incompressible with high probability, as we will now show.

\begin{proposition}
\label{propmainsphere}
Let $u \in \mathbb{S}^{n-1}$ be a uniform random unit vector, and let $0 < \ep < \delta$. Then \[ \mathbf{P}[u \not\in I_n(\delta, \ep)] \le \exp\left\{ 2 \delta (1 + \log 1/\delta) n - \left(\frac{\ep}{\delta} - 1\right)^2n\right\} + \exp\{ - \ep^2 n / 8\}.\]
\end{proposition}
\begin{corollary}
\label{cormainsphere}
Let $u \in \mathbb{S}^{n-1}$ be a uniform random unit vector, and let $0<\delta\leq 1/50$.  Then 
\[\mathbf{P}[u \not\in I_n(\delta, \delta/2)] \le 2\exp\{-\delta^2n/32\}.\]
\end{corollary}
\begin{proof}
Set $\ep = \delta/2$ in Proposition \ref{propmainsphere}. Note that $\varepsilon^2/8 = \delta^2/32$ and  $2\delta(1+\log 1/\delta) - (1/2)^2 < -1/32$ for $0<\delta \leq 1/50$.
\end{proof}

The proof of the Proposition \ref{propmainsphere} consists of two parts. First, we prove a similar proposition where instead of the $u_i$ we have independent Gaussian random variables with the same variance $1/n$. We then use a coupling argument to conclude the desired bound for $u$ drawn uniformly from the unit sphere.

We will need upper and lower tail bounds on the $\chi^2$ distribution.  One can get good enough bounds using the Chernoff method, but rather than develop these from scratch we will cite the following corollary of Lemma 1 from Section 4.1 of \cite{laurent2000adaptive}.
\begin{lemma}\label{lemmachitail}Let $Y$ be distributed as $\chi^2(k)$ for a positive integer $k$.  Then the following upper and lower tail bounds hold for any $t \ge 0$:
\[ \mathbf{P}\left[Y \le k - 2\sqrt{kt}\right] \le e^{-t},\] \[ \mathbf{P}\left[Y \ge k + 2\sqrt{kt} + 2t\right] \le e^{-t}.\] 
\end{lemma}

\begin{proof}[Proof of Proposition \ref{propmainsphere}]
Let $X_1, \dots, X_n$ denote independent Gaussian random variables each with variance $1/n$, and let $X = (X_1, \dots, X_n)$.  If we set $u = X/\Vert X \Vert$, then $u$ is uniformly distributed on the unit sphere; see e.g. \cite{muller1959note}.  

We seek to upper bound the probability of compressibility $\{u \not\in I_n(\delta, \ep)\}$, which is the event that $\sum_{i \in S} u_i^2 < \ep$ for some subset $S$ of coordinates with $|S| \ge \delta n$.  This event is contained in the union of the following two events:

\begin{enumerate}
    \item $E$, the event that $\sum_{i \in S} X_i^2 \le 2\ep$ for some $|S| \ge \delta n$, and 
    \item $F$, the event that $\sum_{i \in S} X_i^2  \ge \ep + \sum_{i \in S} u_i^2$  for some  $|S| \ge \delta n$.
\end{enumerate}
Indeed, if neither of these events hold, then for all $|S| \ge \delta n$ we have \[2\ep < \sum_{i \in S} X_i^2 < \ep + \sum_{i \in S} u_i^2,\] so $u$ is incompressible.

To upper bound the probability of $E$, we use the union bound over all sets of size $k = \lceil n \delta \rceil$: 
\begin{align*}
\mathbf{P}[E] &\le \binom{n}{k} \mathbf{P} \left[ \sum_{i=1}^k X_i^2 \le 2\ep \right] \\ &\le (en/k)^{k} \exp \left\{ - \frac{(k - 2n \ep)^2}{4k} \right\},
\end{align*}
where in the last step we apply the lower tail bound in Lemma \ref{lemmachitail} with $t$ being the solution to $k - 2\sqrt{kt} = 2n\ep$.  To avoid the bookkeeping of ceiling and floor functions we use the extremely crude inequality $n\delta \le k \le 2n\delta$ (valid as long as $\delta n \ge 1$), which will suffice for our purposes:
$$\mathbf{P}[E] \le \exp\left\{2 \delta (1 + \log \delta^{-1} )n -  \left(\frac{\ep}{\delta}- 1\right)^2 n \right\}. $$

We now upper bound the probability of $F$:
\begin{align*}
    \mathbf{P}[F] &= \mathbf{P}\left[\sum_{i \in S} \left(X_i^2 -  \frac{X_i^2}{\Vert X \Vert^2} \right)  \ge \ep  \text{  for some $|S| > \delta n$} \right] \\
    &= \mathbf{P} \left[ \left( 1- \frac{1}{\Vert X \Vert^2} \right) \sum_{i \in S} X_i^2 \ge \ep \text{  for some $|S| > \delta n$} \right] \\
    &\le \mathbf{P} \left[\left(1-  \frac{1}{\Vert X \Vert^2 } \right) \Vert X \Vert^2 \ge \ep \right] \\
    &= \mathbf{P} \left[ \Vert X \Vert^2 \ge 1 + \ep \right]. 
\end{align*}
Since $Y = n \Vert X \Vert^2$ is distributed as $\chi^2(n)$, we may apply the upper tail bound in Lemma \ref{lemmachitail} with $t = n\ep^2/8$ to obtain

\[ \mathbf{P}[F] \le \exp\{-n \ep^2/8\}.\]
To conclude, we have $\mathbf{P}[u \not\in I_n(\delta, \ep)] \le \mathbf{P}[E] + \mathbf{P}[F]$, and substituting the bounds we just derived, we obtain the desired inequality.
\end{proof}

\subsection{\emph{K}-connectedness of the incompressible region}
Having proven that the incompressible region $I_n(\delta, \ep)$, where we have a good local Lipschitz constant, is almost the entire sphere, we now turn to proving that the region is $K$-connected for a small enough $K$.  

One could try to show that any two points in $I_n(\delta, \ep)$ can be connected by a short path contained in $I_n(\delta, \ep)$, but for our purposes it is okay to let the path venture out into the larger region $I_n(4\delta, \ep/\sqrt{2})$. When upgrading to a global Lipschitz constant, we will have to use the slightly worse upper bound for the local Lipschitz constant in this larger region, but this will still be good enough.

\begin{proposition}\label{propkconnected}
$I_n(\delta, \ep)$ is $\sqrt{2/\ep}$-connected in $I_n(4\delta, \ep/\sqrt{2}).$  
\end{proposition}
\begin{proof}
Let $x$ and $y$ be any two endpoints in $I_n(\delta, \ep)$.  The construction will proceed in two steps.  First, we will construct a path from $x$ to $y$ in $\mathbb{R}^n$ consisting of $\lceil \delta^{-1} \rceil$ pairwise orthogonal line segments.  Then we will project this path radially onto the unit sphere and show that the result indeed lies in $I_n(4\delta, \ep/2)$ and has length at most $(2/\sqrt{\ep}) \Vert x - y \Vert$, which is at most $(2/\sqrt{\ep}) d(x,y)$, where $d$ denotes the geodesic distance on $\mathbb{S}^{n-1}$.

Roughly speaking, we will partition the coordinates of $x$ into $1/\delta$ blocks of $\delta n$ coordinates and move the entries of each block linearly from $x$ to $y$ in parallel, one block at a time.

Because basic quantities such as $1/\delta$ and $\delta n$ may not be integers, we will be content to split up $\mathbb{R}^{n}$ as the direct sum $\bigoplus_{i=1}^m \mathbb{R}^{n_i}$, where $\delta n  \le n_i \le 2 \delta n$ for all $i$.\footnote{This is possible as long as $n/2 \ge \delta n \ge 1$, which will be true in our regime.}  Note also that this implies $m \ge \frac{1}{\delta}$.  Similarly, for any vector $z \in \mathbb{R}^n$, we will write $z = \bigoplus_{i=1}^m z^{(i)}$, where $z^{(i)} \in \mathbb{R}^{n_i}.$

Now we may formally define the path $P_i$ to be the line segment 
\[P_i(t) = x^{(1)} \oplus \dots \oplus x^{(i-1)} \oplus \left(t x^{(i)} + (1-t) y^{(i)}\right) \oplus y^{(i+1)} \oplus \dots \oplus y^{(m)}\]
and define $P$ to be the concatenation of the segments $P_1, \dots, P_m$.  The length of $P$ is 
\[ \sum_{i=1}^m \Vert x^{(i)} - y^{(i)} \Vert \le \sqrt{m} \Vert x - y \Vert \le \sqrt{1/\delta} \Vert x - y \Vert,\]
by the Cauchy-Schwarz inequality.  Also, $\Vert P(t) \Vert \ge \sqrt{\ep/2\delta}$, because 
\[\Vert P_i(t) \Vert^2 \ge \sum_{j=1}^{i-1} \Vert x^{(j)}\Vert^2 + \sum_{j=i+1}^{m} \Vert y^{(j)} \Vert^2 \ge (m-1) \ep \ge \frac{\ep}{2\delta},\]
where we use that $x$ and $y$ are $(\delta, \ep)$-incompressible.

Furthermore, note that $P$ lies inside the closed ball of radius $\sqrt{2}$, because for any $i$ and $t$, $$\Vert P_i(t) \Vert^2 \le \sum_{j=1}^{m} \max\{\Vert {x}^{(j)} \Vert, \Vert {y}^{(j)} \Vert\}^2 \le \sum_{j=1}^{m} \left( \Vert {x}^{(j)}  \Vert^2 + \Vert {y}^{(j)}  \Vert ^2 \right) = 2.$$ 

The path $P$ currently does not lie in the unit sphere, so we project it onto the unit sphere along radii to get our final path $P'$. We now show that $P'$ indeed lies in $I_n(4\delta, \ep/\sqrt{2})$.  

At this stage, we will dispense with the direct sum decomposition and use ordinary coordinates $z = (z_1, \dots, z_n)$.  

Consider any set $S$ of at least $4\delta n$ coordinates, and consider any point $P_i(t)$ in our path $P$ (before projection).  The $i$th block of coordinates is in motion, and all of the other coordinates are either frozen at their initial value (from $x$) or their final value (from $y$).

The $i$th block consists of at most $2\delta n$ coordinates.  Besides these, there are at least $4\delta n - 2\delta n = 2\delta n$ remaining coordinates in our set $S$.  At least $\delta n$ of them are from $x$ or at least $\delta n$ of them are from $y$.  By incompressibility of $x$ and $y$, the sum of the squares of these $\delta n$ coordinates is at least $\ep$.

After projecting onto the unit sphere, the sum of the same coordinates is still at least $\ep / \sqrt{2} $, because as we saw, the original path had norm at most $\sqrt{2}$ at every point.

Finally, when projecting onto the unit sphere, the length of the path increases by at most a factor of $1/ \sqrt{\ep/2\delta}$, because as we saw earlier, originally each segment lay outside the smaller sphere of radius $\sqrt{\ep/2\delta}$. The verification is an exercise in plane geometry (using the fact that $\tan \theta > \theta$ for $0 < \theta < \pi/2$) and also follows from the arc length formula $d\textbf{s} = \sqrt{r^2 + (dr/d\theta)^2}\,d\theta \ge r \, d\theta$. 

Thus, finally, we have shown that the path $P'$ is contained in $I_n(4\delta, \ep/\sqrt{2})$ and has length at most
\[\sqrt{1/\delta} \Vert x - y \Vert (1/\sqrt{\ep/2\delta}) = \sqrt{2/\ep} \Vert x - y \Vert.\]
\end{proof}

\section{Concentration of the output}

We now analyze the local Lipschitz constant for the entries $\alpha_i$ and $\beta_i$ of the Jacobi matrix. To simplify notation, in what follows we assume that $\norm{A}_{}=1$ by rescaling $A$. Recall that  this will also rescale the Ritz values and Jacobi coefficients by a factor $1/\norm{A}_{}$.   

By Corollary \ref{propositionhessenberglipschitz}, the function $\alpha_i(u)$ has local Lipschitz constant $2 \cdot 4^i \gamma_{i}(u)$, and $\beta_i(u)$ has local Lipschitz constant $4^{i+1} \gamma_i(u)$. Thus we are naturally led to the question of finding upper bounds for $\gamma_k(u)$.
Recall that $\gamma_k(u)$ is defined as the leading coefficient of the $k$th orthogonal polynomial with respect to the measure $\mu^u = \sum_{i=1}^n u_i^2 \delta_{\lambda_i}$ and that $\pi_k^u$ is the \emph{monic} orthogonal polynomial with respect to the same measure.  

The Equations (\ref{leadingcoeff}) and (\ref{eqmuu}) imply 
\[ \gamma_k(u) = \left( \sum_{i=1}^n u_i^2 \pi_k^u(\lambda_i)^2 \right)^{-\frac{1}{2}}.\]
We seek to upper bound $\gamma_k(u)$ in terms of $u$, so we need to lower bound the quantity 
\[\sum_{i=1}^n u_i^2 \pi_k^u(\lambda_i)^2 = \sum_{i=1}^n  u_i^2 \prod_{j=1}^k |\lambda_i - r_j(u)|^2, \]
where $r_1(u), \dots, r_k(u)$ are the roots of $\pi_k^u(z)$, i.e. the Ritz values.

Now, if it happens to be the case that the $n$ eigenvalues $\lambda_i$ are all clustered very close to the $k$ Ritz values $r_j$, then we won't get a good lower bound.  However, if $k \ll n$ and if the $\lambda_i$ are reasonably spread out, we expect to get a good lower bound for most $i$.  To make this precise, we are led to the notion of equidistribution, which was stated in Section 2.3 and which we restate below.
\\

\begin{definition}(Restatement of Definition \ref{defnequidistribution}). 
Let $\Lambda$ be any finite set of $n$ real numbers.  Let $\delta$ and $\omega$ be positive real numbers and let $j$ be a natural number.  We say that $\Lambda$ is $(\delta, \omega, j)$-\emph{equidistributed} if for any finite set $T$ of at most $j$ real numbers,
\[ \left| \left\{ \lambda \in \Lambda : \frac{1}{|T|} \sum_{t \in T} \log |\lambda - t| \ge \log \omega \right\} \right|  \geq \delta n.\]
\end{definition}

We will show in Section \ref{secequidistribution} that a wide range of spectra are equidistributed. 

Now we apply the definition.  Returning to our effort to upper bound $\gamma_j(u)$, we see that if we assume the spectrum of $A$ is $(\delta, \omega, k)$-equidistributed, then

\[\sum_{i=1}^n  u_i^2 \prod_{j=1}^k |\lambda_i - r_j(u)|^2 \ge \sum_{i \in S} u_i^2 \omega^{2k},\]
where $S$ is some subset of $\{1, \dots, n \}$ of size at least $\delta n$.  However, for an arbitrary unit vector $u$ and an arbitrary subset $S$, we have no lower bound on the sum $\sum_{i \in S} u_i^2$---it could even be zero.  This leads to our definition of incompressibility in Section 3, which is satisfied by $u$ with high probability.

Indeed, if we assume that the unit vector $u$ is $(\delta, \ep)$-incompressible, then the right hand side expression above is greater than $\ep \omega^{2k}$.  Putting together the last few equations, we have $\gamma_k(u) \le (\ep \omega^{2k})^{-1/2}$.  We summarize the result in the following proposition.

\begin{proposition}\label{propgammabound}
Suppose the spectrum of $A$ is $(\delta, \omega, k)$-equidistributed and suppose that $u$ is $(\delta, \ep)$-incompressible for some $\delta, \omega, \ep > 0$ and $k \in \mathbb{N}$. Then
\[\gamma_k(u) \le \frac{1}{\omega^k \sqrt{\ep}}.\]
\end{proposition}

\subsection{Equidistribution}\label{secequidistribution}

In this section we establish sufficient conditions for equidistribution that apply to a wide range of spectra.  First, we present an immediate generalization of the notion of equidistribution which applies to measures $\mu$ instead of finite sets $\Lambda$.  The definitions coincide for finite sets if one identifies $\Lambda$ with the uniform probability distribution on $\Lambda$.
\begin{definition}[Equidistribution for measures]
Let $\mu$ be a probability measure on $\mathbb{R}$.  Let $\delta, \omega > 0$ and $j $ be a natural number.  We say that $\mu$ is $(\delta, \omega, j)$-\emph{equidistribu\-ted} if for any finite set $T$ of at most $j$ real numbers,
\[ \mu \left(\left\{ x \in \mathbb{R} : \frac{1}{|T|} \sum_{t \in T} \log |x - t| \ge \log \omega \right\}\right) \geq \delta.  \]
If a measure is $(\delta,\omega, j)$-equidistributed for every $j \in \mathbb{N}$, we will just say that it is $(\delta, \omega)$-equidistributed. 
\end{definition}

  For absolutely continuous measures, we have the following general equidistribution result.

\begin{proposition}[Absolutely continuous measures are equidistributed] 
\label{propositionAbcContSparse}
Let $\nu$ be a compactly supported probability measure on $\mathbb{R}$ with a nontrivial absolutely continuous part.  Then there exist constants $\delta, \omega > 0$ such that $\nu$ is $(\delta, \omega)$-equidistributed.
\end{proposition}
\begin{proof}
By the assumption, we may write $\nu$ = $\nu_1 + \nu_2$, where $\nu_1$ is absolutely continuous with respect to Lebesgue measure.  By cutting off the portion where the density of $\nu_1$ is greater than some large $M>0$ and assigning that mass to $\nu_2$ instead, we may assume without loss of generality that the density function of $\nu_1$ is bounded.

We now utilize a Markov inequality type argument.  Let $T$ be any set of $j$ real numbers.  Define the logarithmic potential 
$$V(x) =  -\frac{1}{j} \sum_{t \in T} \log |x-t|.$$

Since $\nu_1$ has a bounded density function, $\log |x - t|$ is integrable against $\nu_1$ for all $t$, so the integral $\int_{-\infty}^\infty V_t(x)\,d\nu_1(x)$ is finite for each $t \in T$. Averaging over all $t \in T$, we find that
 
$$\frac{1}{\nu_1(\mathbb{R})} \int_{-\infty}^\infty V(x) d\nu_1(x) \le a$$
for some constant $a < \infty$.  Then
\[a \ge \frac{1}{\nu_1(\mathbb{R})}\int_{-\infty}^\infty V(x) d\nu_1(x) \ge \frac{2a \nu_1(\{x \in \mathbb{R} : V(x) \ge 2a\})}{\nu_1(\mathbb{R})}.\]

Relating this back to the definition of equidistribution, we have
\[ \nu_1\left(\left\{ x \in \mathbb{R} : \frac{1}{|T|} \sum_{t \in T} \log |x - t| \ge -2a \right\}\right) = \nu_1(\{x \in \mathbb{R} : V(x) \le 2a \}) \ge \frac{1}{2} \nu_1(\mathbb{R}).\]
Hence we may take $\delta = \frac{1}{2}\nu_1(\mathbb{R})$ and $\omega = e^{-2a}$.
\end{proof}

Given our framework, it will be useful to have a statement relating the equidistribution of an absolutely continuous measure to a discretization of that measure.  If the two measures are close in Kolmogorov distance, then we can prove such a statement.
\begin{proposition}
\label{propkolconv}
Let $\mu$ and $\nu$ be probability measures.   If $\mu$ is $(\delta, \omega, j)$-equidistri\-buted for some $\delta, \omega > 0$ and $j \in \mathbb{N}$, then $\nu$ is $(\delta - \ep, \omega, j)$-equidistributed, where $\ep =  4j \mathrm{Kol}(\mu, \nu)$.
\end{proposition}
\begin{proof}
Let $T$ be any set of at most $j$ real numbers.  Since $p(x) = \prod_{t \in T} |x - t|$ is the absolute value of a polynomial of degree $j$, each of its level sets is a union of at most $2j$ intervals.  Hence,

$$
|\mu(\{ x\in \mathbb{R} : p(x)\geq \omega^{|T|}\} )- \nu(\{ x \in  \mathbb{R} : p(x)\geq \omega^{|T|}\})| \leq 4j \mathrm{Kol}(\mu, \nu).$$
\end{proof}

Thus, to prove equidistribution for an atomic measure, it suffices to prove equidistribution for a nearby absolutely continuous measure.

The above propositions immediately yield a useful corollary for analyzing the Lanczos procedure in the regime of $O(\log n)$ iterations.
\begin{corollary}\label{corkolconv} 
Let $\mu$ be a compactly supported probability measure with nontrivial absolutely continuous part.  Let $\{\mu_n\}$ be a sequence of probability measures such that $\mathrm{Kol}(\mu_n, \mu) \le \frac{C}{\log n}$ for some $C > 0$. Then for all $n$, for all $j \le \frac{1}{2C} \log n$ we have that $\mu_n$ is $(\delta, \omega, j)$-equidistributed for some $\delta, \omega > 0$.
\end{corollary}

\begin{remark}\label{rem:affine}
If $\mu$ is $(\delta, \omega, j)$-equidistributed and $\nu$ is the pushforward of $\mu$ under the affine map $x \mapsto ax+b$, then $\nu$ is $(\delta, a\omega, j)$-equidistributed.
\end{remark}

We now compute the equidistribution for a few example measures, following the proof of Proposition \ref{propositionAbcContSparse}.  

\begin{example} \label{exampunif} Let $\mu$ denote the uniform measure on $[0,1]$.  Then \[\int V(x) \, d\mu(x)  \le \int -\log\left|x - \frac{1}{2}\right|\,d\mu(x) = 1 + \log 2.\]  Thus, $\mu$ is $(1/2, 4e^{-2})$-equidistributed.
\end{example}

\begin{example} \label{exampsemicirc} Let $\nu$ denote the \emph{semicircle law} $d\nu = \frac{1}{2\pi} \sqrt{(4-x^2)_{+}}\,dx$.  Then \[\int V(x) \, d\nu(x)\le \int -\log\left|x\right|\,d\nu(x) = 1/2.\]  Thus, $\nu$ is $(1/2, e^{-1})$-equidistributed.

\end{example}

With the above the claims made in the examples of Section 2.3 are now trivial. 

\begin{proof}[Proof of Example \ref{exampleparameters} and Example \ref{exampleparameters2}]
It is enough to put together Proposition \ref{propkolconv} and Example \ref{exampunif}. 
\end{proof}

Note that for a given set of points that does not resemble a discretization of an absolutely continuous distribution, it  will still be likely that the equidistribution parameters are well behaved (relative to their scale) provided that the points are somewhat spread out. On the other hand, if the points are clustered in a few small clusters the analysis becomes trivial.   

\begin{observation}
\label{obsboringcase}
Let $\Lambda$ be a set (or multiset) of $n$ points.  Let $a_1 \leq b_1 < a_2\leq b_2 < \cdots < a_m \leq b_m$ be such that  $\Lambda \subset \bigcup_{i=1}^m [a_i, b_i]$. Define $n_i = |\Lambda \cap [a_i, b_i] |$ and let $g$ the minimal gap between clusters, namely, $g = \min_{1\leq i \leq m-1} a_{i+1}-b_i$. Then $\Lambda$ is $(\frac{k_j}{n}, \frac{g}{2}, j)$-distributed, where $k_j = \min_{S} \sum_{i\in S^c} n_i$ and $S$ runs over all subsets of $\{1, \dots, m\}$ of size $j$. 
\end{observation}

\begin{proof}
The proof follows directly from the definition of equidistribution.
\end{proof}

\begin{remark}
A particular case of Observation \ref{obsboringcase} is when $n_i  \geq \lfloor \frac{n}{m} \rfloor$ and $g= a_{i+1}-b_i$  for every $i=1, \dots, m$, which yields Example \ref{exampleparamclusters} above. More generally, if each $n_i$ is roughly $n/m$, then $k_j$ will be roughly $m-j$, and hence the $\delta$ parameter for the equidistribution of $\Lambda$ will  only degrade  when $j \approx m$. In other words, Theorem \ref{thmhessconc} is still strong for matrices whose spectrum consists of small clusters if the number of such clusters exceeds the number of iterations of the Lanczos procedure.  On the other hand, if the number of iterations exceeds the number of clusters it is not hard to show that the Lanczos procedure will output  (with overwhelming probability) at least one Ritz value per cluster.   
\end{remark}

\subsection{Jacobi coefficients}
We now have the necessary tools to prove concentration for the entries of the Jacobi matrix.
\begin{proposition}[Jacobi coefficients are globally Lipschitz]\label{propgloballip}
Suppose the spectrum of $A$ is $(4\delta, \omega, i)$-equidistributed for some $\delta, \omega > 0$ and $i \in \mathbb{N}$. Then for any $0 < \ep < \delta$, functions $\alpha_{i}(u)$ and $\beta_i(u)$ are globally Lipschitz on $I_n(\delta, \varepsilon)$ with constant $L_{i, \ep}\le \frac{4^{i+2}\norm{A}_{{}}^{i+1}}{\omega^i\ep}$.
\end{proposition}

\begin{proof}
Proposition \ref{propositionhessenberglipschitz} says that $\alpha_i(u)$ and $\beta_i(u)$ both have local Lipschitz constant at most $4^{i+1} \Vert A \Vert^{i+1} \gamma_i(u)$ for all $u\in\mathbb{S}^{n-1}$.  Proposition \ref{propgammabound} says that because the spectrum of $A$ is $(4\delta, \omega, i)$-equidistributed, $\gamma_i(u) \le \frac{1}{\omega^i \sqrt{\ep/\sqrt{2}}}$ for all $u \in I_n(4\delta, \ep/\sqrt{2})$.  Combining these, we have that $\alpha_i(u)$ and $\beta_i(u)$ are locally Lipschitz with constant
\[\frac{ 4^{i+1}\norm{A}_{{}}^{i+1}}{\omega^i \sqrt{\ep/\sqrt{2}}}\] 
for all $u \in I_n(4\delta, \ep/\sqrt{2})$.
Proposition \ref{propkconnected} says that $I_n(\delta, \ep)$ is $\sqrt{2/\ep}$-connected in the larger set $I_n(4\delta, \ep/\sqrt{2})$, so Lemma \ref{lemmalocalglobal} implies that $\alpha_i(u)$ and $\beta_i(u)$ are \emph{globally} Lipschitz on $I_n(\delta, \ep)$ with constant
\[ L_{i, \ep} = \frac{\sqrt{2}}{\sqrt{\ep}} \left(\frac{ 4^{i+1}\norm{A}_{{}}^{i+1}}{\omega^i \sqrt{\ep/\sqrt{2}}} \right) \le \frac{4^{i+2}\norm{A}_{{}}^{i+1}}{\omega^i\ep}.  \] 

\end{proof}
We now have the tools to prove our first main theorem, which quantifies the concentration of the Jacobi coefficients around their medians.

\begin{theorem}[Restatement of Theorem \ref{thmhessconc}]
Suppose the spectrum of $A$ is $(\delta, \omega, i)$-equidistributed for some $\delta, \omega > 0$ and $i \in \mathbb{N}$.  Let $\tilde{\alpha}_i$ and $\tilde{\beta}_i$ denote the medians of the Jacobi coefficients $\alpha_i(u)$ and $\beta_i(u)$, respectively. Then for all $t > 0$, the quantities
$\mathbf{P}[|\alpha_{i}(u) - \Tilde{\alpha_{i}}| > t\Vert A \Vert]$ and $\mathbf{P}[|\beta_{i}(u) - \Tilde{\beta_{i}}| > t\Vert A \Vert]]$ are both bounded above by 
\begin{equation}
2\exp\left\{-\frac{\min\{\delta, 1/50\}^2}{32}n\right\} + 2 \exp \left\{- \frac{1}{64}\left(\frac{\omega}{4 \Vert A \Vert} \right)^{2i}\delta^2 t^2 n\right\}. 
\end{equation}
\end{theorem}

\begin{proof}
The local L\'evy lemma (Lemma \ref{lemmalevy}) yields that $\mathbf{P}[  |\alpha_i(u) - \tilde{\alpha_i}| > t\Vert A \Vert_{{}}]$ and $\mathbf{P}[  |\beta_i(u) - \tilde{\beta_i}| > t\Vert A \Vert_{{}}]$ are both at most
\[\mathbf{P}[u \not\in I_n(\delta
, \ep)] + 2 \exp\{- 4 n t^2 \Vert A \Vert^2 / L_{i, \ep}^2\}, \]
where $L_{i, \ep}$ is the global Lipschitz constant on $I_n(\delta, \ep)$ obtained in Proposition \ref{propgloballip}.  Note that if $\delta > 1/50$, then $A$ is still $(1/50, \omega, i)$-equidistributed, so we may set $\ep = \delta/7$ and apply Corollary \ref{cormainsphere} to bound $\mathbf{P}[u \not\in I_n(\delta
, \ep)]$.  We obtain the upper bound
\[ 2\exp\left\{-\frac{\min\{\delta, 1/50\}^2}{32}n\right\} + 2 \exp \left\{ \frac{-4nt^2 \Vert A \Vert^2 \omega^{2i} (\delta/2)^2 }{4^{2i+4} \Vert A \Vert^{2i+2}} \right\}\]\[\le 2\exp\left\{-\frac{\min\{\delta, 1/50\}^2}{32}n\right\} + 2 \exp \left\{- \frac{1}{64}\left(\frac{\omega}{4 \Vert A \Vert} \right)^{2i}\delta^2 t^2 n\right\}\]
as desired.
\end{proof}

Combining the previous theorem with Corollary \ref{corkolconv} we get convergence in probability of the Jacobi matrices in the regime $k = O(\log n)$.

\begin{proposition}\label{prophesslogconvprob}
Let the spectra $\mu_n$ of $A_n$ converge to the spectrum $\mu$ of $A$ in Kolmogorov distance with rate $O(1/\log n)$.  Suppose $\mu$ has a nontrivial absolutely continuous part.  Then there exists $c_2 > 0$ and a sequence $k_n \ge c_2 \log n$ such that the Jacobi matrices $J_{k_n}$ output by the Lanczos algorithm after $k_n$ iterations converge to entrywise in probability to deterministic constants.
\end{proposition}
\begin{proof}
By Corollary \ref{corkolconv}, we have that $\mu_n$ is $(\delta, \omega, k)$-equidistributed for all $k \le c_1 \log n$.  Picking $c_2 < c_1$ and applying Theorem \ref{thmhessconc}, for $i \le c_2 \log n$ this yields the bound
\begin{align*} \mathbf{P}[|\alpha_{i} - \Tilde{\alpha_{i}}| > t] 
&\le \exp\{-\delta^2n/32\} + 2 \exp\left\{- \frac{4}{4^3} (\omega/4)^{2 c_2 \log n} n t^2\right\} \\
 &= \exp\{-\delta^2n/32\} + 2 \exp\left\{- \frac{4}{4^3} n^{2 c_2 \log(\omega/4) + 1} t^2\right\} 
\end{align*}
so as long as $ 2 c_2 \log(\omega/4) + 1 > 0$, we have convergence in probability of the Jacobi coefficients as $n \to \infty$.  But this is certainly true for small enough $c_1$.  The $\beta_i$ have the same bound as the $\alpha_i$, so we are done.
\end{proof}

As mentioned in the introduction, convergence for \emph{fixed} $k$ to the infinite Jacobi matrix $J$ of $\mu$ for deterministic $\mu_n$ (with no hypothesis on the rate of convergence of $\mu_n$) is proven in \cite[Theorem 4]{gautschi1968construction}.  In Proposition \ref{prophesslogconvprob} we leave it open to prove that the limit is actually $J$ (see Section \ref{sec:conclusion}), but if we reduce the number of iterations from $k = O(\log n)$ to $k = O(\sqrt{\log n})$, we can indeed prove that the limit is $J$.  This is the content of Theorem \ref{Thmsqrtiter}, proven in Section 5.

\subsection{Ritz values} 
 \label{subsec:RitzVals}
Theorem \ref{thmhessconc} yields  concentration of the entries of the random matrix $J_k(u)$. In general, controlling the entries of a random matrix does not yield control over its random eigenvalues, but, since $J_k(u)$ is Hermitian we know that its spectrum is stable with respect to small perturbations of the entries. More precisely, we will use the well known Weyl's inequality---see \cite[Theorem 4.3.1]{horn2012matrix} for a reference.

\begin{lemma}[Weyl]
\label{lemmaWielandt-Hoffman}
For every matrix $X$, let $\lambda_1(X)\geq \dots \geq \lambda_n(X)$ denote the eigenvalues of $X$. 
If $A$ and $B$ are $n \times n$ Hermitian matrices, then for all $1 \le i \le n$ we have
$$|\lambda_i(A+B)-\lambda_i(A) |\leq \norm{B}_{{}}.$$
\end{lemma}

Following the notation in Theorem \ref{thmhessconc}, let $\tilde{J}_k$ be the $k\times k$ Jacobi matrix with entries $\tilde{\alpha}_i$ and $\tilde{\beta}_i$, and denote the eigenvalues of $\tilde{J}_k$ by $\tilde{r}_1\geq \dots \geq \tilde{r}_k$. 

\begin{proposition}[Concentration of the Ritz values]
\label{propConcRitzVals}
Assume that the spectrum of $A$ is $(\delta, \omega, k)$-equidistributed for some $\delta, \omega > 0$ and $k \in \mathbb{N}$. With the notation described above, let $\vec{r} = (\tilde{r}_1, \dots, \tilde{r}_k) $ and let $\vec{r}(u) = (r_1(u),$ $\dots, r_k(u))$ be the vector of Ritz values after $k$ iterations. Then the probability $\mathbf{P}[\norm{\vec{r}(u)-\vec{r}}_\infty \geq t\Vert A \Vert_{{}}]$ is  bounded above by
\begin{align*}
  4k \left[ \exp\left\{-\frac{\min\{\delta,  1/50\}^2}{32}n\right\} + \exp\left\{-\frac{1 }{192} \left(\frac{\omega }{4 \norm{A}_{{}}}\right)^{2k}  \delta^2 t^2n\right\} \right]. 
\end{align*}

\end{proposition}

\begin{proof}
  Since $\tilde{J}_k$ and $J_k(u)$ are tridiagonal matrices, we may split $J_k - \tilde{J}_k$ into the sum of three matrices consisting of the diagonal, the subdiagonal, and the superdiagonal and then use the triangle inequality to obtain
\begin{equation}
\label{eqq}
\norm{J_k(u)- \tilde{J}_k}_{{}} \leq  \max_{0\leq i \leq k-1} \{|\alpha_i(u)-\tilde{\alpha}_i|\}+ 2\max_{0\leq i \leq k-2}\{|\beta_i(u)- \tilde{\beta}_i|\}.     
\end{equation}
Hence, we deduce that
\begin{align*}
    \mathbf{P}[\norm{\vec{r}(u)-\vec{r}}_{\infty}\geq t]
    &\leq \mathbf{P}[\norm{J_k(u)-\tilde{J}_k}_{{}} \geq t ] \\
    &\leq \mathbf{P}\left[\max_{0\leq i \leq k-1}\{|\alpha_i(u)-\tilde{\alpha}_i|\}+ 2\max_{0\leq i \leq k-2} \{|\beta_i(u)-\tilde{\beta}_i|\}\geq t \right],
\end{align*}
where the first inequality follows from Lemma \ref{lemmaWielandt-Hoffman} and the second inequality from (\ref{eqq}). Now observe that the  event $\{\max_{0\leq i \leq k-1}\{|\alpha_i(u)-\tilde{\alpha}_i|\}+ 2\max_{0\leq i \leq k-2} \{|\beta_i(u)-\tilde{\beta}_i|\}\geq t\}$ is contained in the event
\begin{align*}
 \left\{\max_{0\leq i \leq k-1}\{|\alpha_i(u)-\tilde{\alpha}_i|\}\geq \frac{t}{3}\right\} \cup \left\{\max_{0\leq i \leq k-2}\{|\beta_i(u)-\tilde{\beta}_i|\}\geq \frac{t}{3}\right\},  
\end{align*}
which in turn is contained in the event $$\bigcup_{i=1}^k \left\{|\alpha_i(u)-\tilde{\alpha}_i|\geq \frac{t}{3} \right\} \bigcup \left\{|\beta_i(u)-\tilde{\beta}_i|\geq \frac{t}{3} \right\}.$$
Using a union bound and applying Theorem \ref{thmhessconc}, we obtain the desired result.
\end{proof}

\subsection{Ritz Vectors} Here we will use the same notation  as in Section \ref{subsec:RitzVals}.  Let $\tilde{w}_i$ be the eigenvector of $\tilde{J}_k$ corresponding to $\tilde{r}_i$ and let $w_i(u)$ be the eigenvector of $J_k(u)$ corresponding to $r_i(u)$. We will use the fact that $J_k(u)$ concentrates around $\tilde{J}_k$, together with the Davis-Kahan theorem \cite{davis1969some} to establish the concentration of the vectors $w_i(u)$. 

\begin{theorem}[Davis-Kahan]
\label{thmDavis-Kahan}
Here we use the notation of Lemma \ref{lemmaWielandt-Hoffman}. Fix $i\in \{1, \dots, n\}$ and assume that $\lambda_i(A)$ has multiplicity 1. Define
$$\varepsilon = \min_{j: j\neq i}|\lambda_{i}(A)-\lambda_j(A)|  ,$$
and let $\theta\in [0, \pi/2]$ denote the angle between the $i$-th eigenvectors of $A$ and $A+B$. Then 
$$\sin \theta \leq \frac{2 \norm{B}}{\varepsilon}.$$
\end{theorem}

Under the assumption that $\tilde{r}_i(u)$ is not close to the other Ritz values, we get the following result. 

\begin{proposition}[Concentration of the Ritz vectors]
\label{propConcRitzeigenvectors} 
Assume that the spectrum of $A$ is $(\delta, \omega, k)$-equidistributed  for some $\delta, \omega > 0$ and $k \in \mathbb{N}$ and fix some $i \in \mathbb{N}$ with $1\leq i\leq k$. With the notation described above, let $\theta \in [0, \pi/2]$ be the angle between $w_i(u)$ and $\tilde{w}_i$ and let $\varepsilon = \min_{j: j\neq i} |\tilde{r}_i-\tilde{r}_j|$. Then for any $0\leq c < 1/2$, the probability $\mathbf{P}\left[\sin \theta \geq 2\Vert A \Vert/\varepsilon n^c\right]$ is bounded above by
\begin{align*}
 4k \left[ \exp\left\{-\frac{\min\{\delta,  1/50\}^2}{32}n\right\} + \exp\left\{-\frac{1 }{192} \left(\frac{\omega }{4 \norm{A}_{{}}}\right)^{2k}  \delta^2 n^{1-2c}\right\} \right]. 
\end{align*}
\end{proposition} 

\noindent Note. The same result holds for the Ritz vectors, since these are obtained by applying an isometry to the $w_i(u)$.

\begin{proof}
From Theorem \ref{thmDavis-Kahan} we have that $$\sin \theta \leq \frac{2\norm{\tilde{J}_k(u)-\tilde{J}_k(u)}}{\varepsilon}$$ and hence 
\begin{align*}
\mathbf{P}[\sin \theta \geq t] & \leq \mathbf{P}[\norm{J_k(u)-\tilde{J}_k} \geq t] \\ & \leq \mathbf{P}\left[\max_{0\leq i \leq k-1}\{|\alpha_i(u)-\tilde{\alpha}_i|\}+ 2\max_{0\leq i \leq k-2} \{|\beta_i(u)-\tilde{\beta}_i|\}\geq t \right],
\end{align*}
where the latter inequality was established in the proof of Proposition \ref{propConcRitzVals}. Using the bound obtained in the aforementioned proof and substituting $t = \frac{2}{\varepsilon n^c}$ we obtain the desired result.   
\end{proof}

\section{Proofs of Proposition \ref{prop:outliers} and Theorem \ref{Thmsqrtiter}}

\subsection{Proof of Proposition \ref{prop:outliers}}
We now prove our theorem about the Lanczos algorithm missing outliers in the spectrum.
\begin{proof}[Proof of Proposition \ref{prop:outliers}]
By Proposition \ref{propkolconv}, we have that $\mu_n$ is $(\delta, \omega, j)$-equi\-distributed for some $\delta, \omega > 0$ and all $j < c \log n$.  Suppose $u \in I_n(\delta, \ep)$, which happens with overwhelming probability by Proposition \ref{propmainsphere}.  Then by Proposition \ref{propgammabound}, we have an upper bound on the leading coefficient of the $j$th orthogonal polynomial: $\gamma_j(u) \le \frac{1}{\omega^j \sqrt{\ep}}.$  Equivalently, this is a lower bound on the $L^2$ norm of the $j$th \emph{monic} orthogonal polynomial: $\Vert\pi_j^u\Vert_{L^2(\mu^u)}  \ge \omega^j \sqrt{\ep}.$  As mentioned in the preliminaries in Section 2, it is a classical fact that the monic orthogonal polynomial of any given degree has minimal $L^2$ norm over all monic polynomials of that degree.  Thus, we in fact have 
\begin{equation}\label{eqnmonicbound}
    \int q(x)^2\,d\mu^u(x) \ge \ep \omega^{2j}
\end{equation}
for all monic polynomials $q$ of degree $j$, with equality when $q(x)$ is the $k$th orthogonal polynomial $p_k^u(x)$.

For all unit vectors $u$, let $\rho(u)$ denote the top Ritz value, i.e. the maximum root of $p_k^u(x)$.  We wish to show that $\rho(u) < R + \kappa$ with high probability.

Take $p_k^u(x)$ and replace its top root by $t$ to form the monic polynomial $P_t$.  By the first-order condition for the  variational characterization of $p_k^u$ mentioned above, to show $\rho(u) \le R+\kappa$ it suffices to show that $\Vert P_t \Vert_{L^2(\mu^u)}$ is strictly increasing in $t$ for $t > R+\kappa$.  We have

\[\Vert P_t \Vert_{L^2(\mu^u)}^2 = \int \left( \frac{\pi_k^u(x)}{x-\rho(u)} (x-t) \right)^2 d \mu^u (x) = \sum_{i=1}^{k} u_i^2(\lambda_i - t)^2 \prod_{j=2}^{k}(\lambda_i - r_j)^2,\]
where we let $r_2, \dots, r_{k}$ denote the roots of $p_k^u(x)$ besides the maximum root $\rho(u)$, and we omit the argument $u$ for brevity.  We calculate the derivative

\[ \frac{d}{dt} \Vert P_t \Vert_{L^2(\mu^u)}^2 = -2 \sum_{i=1}^m u_i^2 (\lambda_i- t) \prod_{j=1}^{k-1}(\lambda_i - r_j)^2 
- 2 \sum_{i=m+1}^n u_i^2 (\lambda_i - t) \prod_{j=2}^{k} (\lambda_i - r_j)^2. \]

We wish to show that this quantity is positive whenever $t \ge R+\kappa$.  We have assumed that there are only $m$ outliers, so assume $\lambda_i \le R$ for all $i > m$.  Then $t - \lambda_i \ge \kappa$ for every $m < i \le n$. 

Thus,
\begin{align*} \frac{d}{dt} \Vert P_t \Vert_{L^2(\mu^u)}^2 &\ge -2 \sum_{i=1}^m u_i^2 (\lambda_i- t) \prod_{j=1}^{k-1}(\lambda_i - r_j)^2 
+ 2 \sum_{i=m+1}^n u_i^2 \kappa \prod_{j=2}^{k} (\lambda_i - r_j)^2 \\
&= -2 \sum_{i=1}^m u_i^2 (\lambda_i- t) \prod_{j=2}^{k}(\lambda_i - r_j)^2 \\ &\qquad + \left[
 2 \kappa \int \left(\frac{p_k^u(x)}{x-\rho(u)}\right)^2 \,d\mu^u(x) 
- 2 \sum_{i=1}^m u_i^2 \kappa \prod_{j=2}^{k} (\lambda_i - r_j)^2 \right]\\
&\ge -2 \sum_{i=1}^m u_i^2 (\lambda_i- t) \prod_{j=2}^{k}(\lambda_i - r_j)^2 + 2 \kappa \ep \omega^{2 (k-1)} - 2 \sum_{i=1}^m u_i^2 \kappa \prod_{j=2}^{k} (\lambda_i - r_j)^2,
\end{align*}
where in the last step we used the inequality (\ref{eqnmonicbound}) on the degree $k-1$ polynomial $p_k^u(x)/(x - \rho(u))$.
Simplifying, we have 

\[ \frac{d}{dt} \Vert P_t \Vert_{L^2(\mu^u)}^2 \ge 2 \kappa \ep \omega^{2 (k-1)} -2  \sum_{i=1}^m u_i^2 (\lambda_i + \kappa - t) \prod_{j=2}^{k} (\lambda_i - r_j)^2.\]
  By uniform boundedness of the spectra, there exists $M$ large such that $\lambda_i - r_j \le M$ for all $1\le i \le m$. Let $g$ be the maximum of the outlier gaps $\lambda_i - R$ over all $1 \le i \le m$.  Recall that $t \ge R + \kappa$, so $\lambda_i + \kappa - t \le \lambda_i - R \le g$ for all $1 \le i \le m$. Finally, we have with overwhelming probability $\sum_{i=1}^m u_i^2 < n^{-c}$ for any positive $c < 1/2$; we will defer the proof to Lemma \ref{lemma:supercompressibility} below.  Putting this all together, we have
\[ \frac{d}{dt} \Vert P_t \Vert_{L^2(\mu^u)}^2 \ge 2 \kappa \ep \omega^{2  k - 2}  - 2n^{-c} M^{2k - 2} mg
.\]
This quantity is strictly positive when
\[\log \kappa \ep + (2k-2) \log \omega > -c \log n + (2k-2) \log M + \log mg.\]
  Rearranging, we get
\[(2k-2) \log (\omega/M) > -c \log n + \log mg - \log \kappa \ep \]
for $n$ large.  Note that $\omega < M$, because $\omega$ is a lower bound on geometric means of distances that are all less than $M$.  In conclusion, with high probability, $\frac{d}{dt} \Vert P_t \Vert_{L^2(\mu^u)}^2 > 0$ for all $t > R+\kappa$ when
\begin{equation}\label{eqn:nonasy}
    2k-2 < \frac{1}{\log\frac{M}{\omega}} \left(c \log n + \log \frac{\kappa \ep}{mg} \right).
\end{equation}

For $n$ large, we may absorb the constants $m, g, \kappa, \ep, \omega$ (which do not depend on $n$) into a single constant $c'>0$, and we get the desired $k \le c' \log n$.

\end{proof}
\begin{remark}
We have focused on the right side of the spectrum for ease of exposition.  Similar results hold for outliers on both sides.
\end{remark}
\begin{remark}
There are several parameters that can be tuned in the above proof.  For example, one could envision a situation in which $\kappa$ converges to zero as $n\to \infty$, at the expense of some other parameter.
\end{remark}
\begin{lemma}\label{lemma:supercompressibility}
Let $0 < c < 1/2$ and suppose  $m \le n^\alpha$, where $\alpha < 1 - c$.  Then $\sum_{i=1}^m u_i^2 < n^{-c}$ with overwhelming probability.  To be precise, 
\[\mathbf{P}\left[\sum_{i=1}^m u_i^2 \ge n^{-c}\right] \le \exp \left\{ - \frac{1}{16} \left(4n^\alpha - 4\sqrt{2} n^{\frac{1}{2} - \frac{c}{2} + \frac{\alpha}{2}} + 2n^{1-c}\right) \right\} + \exp\left\{-\frac{1}{16}n^{1-2c}\right\}.\]
\end{lemma}
\begin{proof}
We proceed just as in the proof of Proposition \ref{propmainsphere}.  Define $X_i$ as in that proof.  Then 
\[ \mathbf{P}\left[\sum_{i=1}^m u_i^2 > n^{-c}\right] \le \mathbf{P}\left[\sum_{i=1}^m X_i^2 > \frac{1}{2} n^{-c}\right] + \mathbf{P}\left[\sum_{i=1}^m X_i^2 < -\frac{1}{2} n^{-c} + \sum_{i=1}^m u_i^2\right]. \]
Using Lemma \ref{lemmachitail}, we solve for the parameter $\sqrt{t}=\frac{-2\sqrt{m} + \sqrt{2} n^{\frac{1}{2} - \frac{c}{2}}}{4}$ (which requires $\alpha < 1-c$) and then we get
\begin{align*}
 \mathbf{P}\left[\sum_{i=1}^m X_i^2 > \frac{1}{2} n^{-c}\right] & \leq \exp\left\{ - \left( \frac{-2\sqrt{m} + \sqrt{2} n^{\frac{1}{2} - \frac{c}{2}}}{4}\right)^2\right\} \\ &   = \exp \left\{ - \frac{1}{16} \left(4n^\alpha - 4\sqrt{2} n^{\frac{1}{2} - \frac{c}{2} + \frac{\alpha}{2}} + 2n^{1-c}\right) \right\},
\end{align*}

which is an overwhelmingly small probability because $\frac{1}{2} - \frac{c}{2} + \frac{\alpha}{2} < 1-c$ when $\alpha < 1-c$.

Now following the same coupling argument in the proof of Proposition \ref{propmainsphere} and using Lemma \ref{lemmachitail} again, we get
\[ \mathbf{P}\left[\sum_{i=1}^m X_i^2 < -\frac{1}{2} n^{-c} + \sum_{i=1}^m u_i^2\right] \le \exp\left\{-\frac{1}{16}n^{1-2c} \right\}.\]
\end{proof}

\begin{proof}[Proof of Theorem \ref{thm:nonasylocation}]
From the proof of Proposition \ref{prop:outliers}, setting $\ep = \delta/2$ we have that the Ritz values are contained in the desired interval for
\[k \le \frac{1}{2 \log \frac{M}{\omega}}\left(c \log n + \log \frac{\kappa \delta}{2 m g}\right) \]
as long as $k \le j$, $u \in I_n(\delta, \delta/2)$ and $\sum_{i=1}^m u_i^2 > n^{-c}$.  Applying Corollary \ref{cormainsphere}, the probability that $u$ violates either condition is at most

\begin{align*}& \textbf{P}[u \not\in I_n(\delta, \delta/2)] + \mathbf{P}\left[\sum_{i=1}^m u_i^2 > n^{-c}\right]  \\ & \le 2\exp\left\{-\frac{\min\{\delta,  1/50\}^2}{32} n\right\} + \mathbf{P}\left[\sum_{i=1}^m u_i^2 > n^{-c}\right]\\
&\le 2\exp\left\{-\frac{\min\{\delta,  1/50\}^2}{32} n\right\} + 2\exp \left\{ -\frac{1}{16} n^{1-2c} \right\},
\end{align*}
where in the last step, we apply Lemma \ref{lemma:supercompressibility} and note that for $n \ge e^\frac{1}{1 - c - \alpha}$ we have $ 4\sqrt{2} n^\frac{1-c+\alpha}{2} \le n^{1-c}$.

\end{proof}

\subsection{Proof of Theorem \ref{Thmsqrtiter}}
 For $C >0$ let $\mathcal{P}_C$ denote  the space of Borel probability measures supported on $[-C, C]$. In order to prove Theorem \ref{Thmsqrtiter} we will show that the Jacobi coefficients  of a measure are locally Lipschitz quantities on the space $\mathcal{P}_C$ equipped with the Kolmogorov metric. Note that in Section 3 similar results were obtained in the case in which the space of measures in consideration is restricted to atomic measures supported on $n$ fixed points, namely, the eigenvalues of $A_n$. Since $\mathcal{P}_C$ is a much larger and complicated space we are not able to obtain  results as strong as in Proposition \ref{propositionhessenberglipschitz}. It remains an open question if a better rate can be achieved at this level of generality; see the concluding remarks for some natural directions to pursue.

We will use the following well known result which, for convenience of the reader, we restate as it appears in Lemma 1.1 in \cite{godsil2017algebraic}. 

\begin{lemma}
\label{lemmaGodsil}
Let $A$ and $B$ be two $k\times k$ matrices. Then $\det(A+B)$ is equal to the sum of the determinants of the $2^k$ matrices obtained by replacing each subset of the columns of $A$ by the corresponding subset of the columns of $B$. 
\end{lemma}
\begin{proof}
The result follows directly from the fact that the determinant is multilinear in the columns of the matrix.
\end{proof}
\begin{lemma}\label{inequalityCorollary}
Let $A$ and $B$ be two $k\times k$ matrices. For $1\leq i \leq k$, let $A^{(i)}$ and $B^{(i)}$ be the $i$th columns of $A$ and $B$, respectively. Let $C, \varepsilon >0$ and assume that 
\begin{equation}\label{assncorollary}
\norm{A^{(i)}-B^{(i)}}_2 \leq \varepsilon \hspace{.3cm} \text{and} \hspace{.3cm}\max\{\norm{A^{(i)}}_2, \norm{B^{(i)}}_2\} \leq C. 
\end{equation}
Then 
$$|\det(A)-\det(B)| \leq \varepsilon k (C+\varepsilon)^{k-1}.$$
\end{lemma}

\begin{proof}
By the assumption in (\ref{assncorollary}) we can write $B = A+E$, where $E$ is a matrix with columns of norm less than or equal to $\varepsilon$. Then, using Lemma \ref{lemmaGodsil}, the inequalities in (\ref{assncorollary}),  and the fact that the determinant of a matrix is bounded by the product of the Euclidean norms of its columns, we obtain
$$|\det(A+E)-\det(A)| \leq \sum_{k=1}^n  {\binom{n}{k}} C^{n-k}\varepsilon^k = (C+\varepsilon)^k- C^k \leq \varepsilon k (C+\varepsilon)^{k-1},$$
where the last inequality follows from the mean value theorem. 
\end{proof}

We now argue that the moments of a measure are Lipschitz quantities in $\mathcal{P}_C$, where the constant is exponential in the order of the moment. With this end fix a Borel measure $\mu$ on $\mathbb{R}$ and denote 
$$m_k(\mu) = \int_\mathbb{R} x^k d \mu(x).  $$
A standard application of Fubini's theorem yields that if $\mu$ is a finite positive Borel measure supported in $[0, \infty)$, then 
\begin{equation}
\label{eqmoments}
m_k(\mu) = k \int_0^\infty  x^{k-1} \mu(x, \infty) d x.
\end{equation}
This identity is enough to obtain the following bound. 

\begin{lemma}
\label{lemmamomentkolmogorov}
Let $\mu, \nu \in \mathcal{P}_C $ and $k>0$, then $|m_k(\mu)-m_k(\nu)| \leq  2 C^k \mathrm{Kol}(\mu, \nu)$. 
\end{lemma}

\begin{proof}
Start by decomposing $\mu$ into $\mu_+$ and $\mu_-$ as follows: 
$$\mu_+(A) = \mu(A \cap [0, \infty)), \hspace{.2cm} \mu_-(A) = \mu(-A\cap (-\infty, 0) ) \hspace{.3cm} \forall A \in \mathcal{B}(\mathbb{R}). $$
Hence $\mu(A) = \mu_+(A) + \mu_-(-A)$.  Define $\nu_+$ and $\nu_-$ analogously. Note that these new measures are supported on $[0, \infty)$. 

Observe that $m_k(\mu) = m_k(\mu_+)+ (-1)^km_k(\mu_-)$ and that the analogous formula holds for  $m_k(\nu)$. Hence 
$$|m_k(\mu)-m_k(\nu)| \leq |m_k(\mu_+)-m_k(\nu_+)|+ |m_k(\mu_-)-m_k(\nu_-)|. $$

Now, for $t\geq 0$ define $F_{\mu_+}(t) =  \mu_+(t, \infty)$ and $F_{\nu_+}(t) = \nu_+(t, \infty)$. By definition of Kolmogorov distance we have that
$$|F_{\mu_+}(t)-F_{\nu_+}(t)|\leq  \mathrm{Kol}(\mu, \nu).$$
On the other hand, by  (\ref{eqmoments}) we have that 
\begin{align*}
 |m_k(\mu_+)-m_k(\nu_+) | &\leq k  \int_0^\infty x^{k-1} |F_{\mu_+}(x)-F_{\nu_+}(x)|  dx \\ & \leq k \mathrm{Kol}(\mu, \nu) \int_0^C x^{k-1} dx \\ & = C^k \mathrm{Kol}(\mu, \nu). 
\end{align*}

In the exact same way we can bound $|m_k(\mu_-)-m_k(\nu_-)|$ to conclude the proof. 
\end{proof}

Given $\mu \in \mathcal{P}_C$  we denote the $(k+1)\times (k+1)$ Hankel matrix of $\mu$ by $M_k(\mu)$ and  define $D_k(\mu) = \det M_k(\mu)$. We will denote the Jacobi coefficients of $\mu$ by $\alpha_i^\mu$ and $\beta_i^\mu$.  For the proof of the following results, many of the facts stated in Section 2.1 will be used.

\begin{proposition}
\label{propositionloclipbeta}
Let $\mu, \nu \in \mathcal{P}_C$ and let $ s_k>0$ be constants satisfying 
\[\min\{D_j(\mu), D_j(\nu)\} \geq s_k\] for $j=1, \dots, k$. Then  
$$|\beta_k^\mu- \beta_k^\nu| \leq  \frac{ \exp\{g  k^2\} \mathrm{Kol}(\mu, \nu)}{s_k^2} $$
 for some $g>0$ dependent of $\mu$ and $\nu$ but independent  of $k$. 
\end{proposition}

\begin{proof} To shorten notation let $x_j = D_j(\mu)$ and $y_j= D_j(\nu)$. Without loss of generality $ C >1$. A direct application of Lemma  \ref{lemmamomentkolmogorov} yields a rough bound between the distance in the Euclidean norm of the corresponding columns of the matrices $M_j(\mu)$ and $M_j(\nu)$. Namely, the columns are at distance less than $ \sqrt{j+1}C^{2j-1} \mathrm{Kol}(\mu, \nu)$. The same reasoning yields that the norm of any column in $M_j(\mu)$ or $M_j(\nu)$ is bounded by $\sqrt{j+1}C^{2j-1}$.  Hence, using  Lemma \ref{inequalityCorollary} we get 

$$|x_j-y_j| \leq (\sqrt{j+1})^{{j+1}} j (C^{(2j-1)} + \ep)^{j+1} \mathrm{Kol}(\mu, \nu) \leq \exp\{g j^2\} \mathrm{Kol}(\mu, \nu)$$
for some $g >0$ independent of $k$. 

In what follows we will bound two other terms whose logarithm is also $O(k^2)$.  The implied constants depend only on $\mu$ and $\nu$, so we can modify $g$ to be big enough for the following inequalities to hold as well.  By the first expression in (\ref{eqorthpolidentities}) we have that
\begin{align}
\label{ineqbetas}
|\beta_k^\mu- \beta_k^\nu| &= \left| \frac{\sqrt{x_{k-1}x_{k+1}}}{x_k} - \frac{\sqrt{y_{k-1}y_{k+1}}}{y_k} \right|\nonumber \\ &\leq \frac{1}{x_k} |\sqrt{x_{k-1}x_{k+1}}- \sqrt{y_{k-1}y_{k+1}}| +\sqrt{y_{k-1}y_{k+1}} \left| \frac{1}{x_k}- \frac{1}{y_k} \right|.\end{align}
To bound the first term on the right-hand side of the above inequality we see that
\begin{align*}
|\sqrt{x_{k-1}x_{k+1}}- \sqrt{y_{k-1} y_{k+1}}| & =  \frac{|x_{k-1}x_{k+1}- y_{k-1} y_{k+1}|}{\sqrt{x_{k-1}x_{k+1}}+ \sqrt{y_{k-1} y_{k+1}}}\hspace{.4cm} \text{and} \\   |x_{k-1}x_{k+1}- y_{k-1}y_{k+1}| & \leq x_{k-1} |x_{k+1}-y_{k+1}| +y_{k+1} |x_{k-1}-y_{k-1}| \\ & \leq \exp\{a k^2\}\mathrm{Kol}(\mu, \nu),
\end{align*}
which yields 
\begin{equation}
\label{ineq}
 \frac{1}{x_k}|\sqrt{x_{k-1}x_{k+1}}- \sqrt{y_{k-1}y_{k+1}}|\leq \frac{\exp\{g k^2\}\mathrm{Kol}(\mu, \nu)}{2s_k^2}.  
\end{equation}
On the other hand,
\begin{equation}
\label{ineqy's}
\sqrt{y_{k-1}y_{k+1}} \left| \frac{1}{x_k}- \frac{1}{y_k} \right| = \sqrt{y_{k-1}y_{k+1}} \frac{|x_k-y_k|}{x_k y_k} \leq \frac{ \exp\{g k^2\}\mathrm{Kol}(\mu, \nu)}{2 s_k^2}.
\end{equation}
The result then follows from combining the previous inequalities (\ref{ineqbetas}), (\ref{ineq}), and (\ref{ineqy's}).  
\end{proof}

\begin{remark}
The constants $s_k$ have already been studied with sophisticated techniques for some families of measures; see \cite{Szego} for an example. However, using results only from Section 4 it will be easy to show that for measures with an absolutely continuous part we have $|\log(s_k)| = O(k^2)$, where the implied constant depends only on $\mu$, which is enough for the proof of Theorem \ref{Thmsqrtiter}.  
\end{remark}

In a similar fashion we can show that the coefficients of $p_k^\mu(x)$ are locally Lipschitz. 

\begin{proposition}
\label{proplocalLipscoefficients}
Fix a positive integer $k$. Let $\mu, \nu$ and $s_k$ be  as in Proposition \ref{propositionloclipbeta}. Denote the coefficients of $x^i$ in $p_k^\mu(x)$ and $p_k^\nu(x)$ by $a_i^\mu$ and $a_i^\nu$ respectively. Then
$$|a_i^\mu- a_i^\nu| \leq  \left(\frac{2}{s_k}+ \frac{1}{s_k^2}\right) \mathrm{Kol}(\mu, \nu) \exp\{g k^2\}$$
for some $g>0$ dependent on $\mu$ and $\nu$ but independent  of $k$. 
\end{proposition}

\begin{proof}
For $1\leq i \leq k$ let $M_k^{(i)}(\mu)$ be the matrix obtained by removing the $k$th row and $i$th column  of $M_k(\mu)$ and let $d_i(\mu) = \det(M_k^{(i)}(\mu))$. From identity (\ref{eqdeterminantalid}) we have 
$$a_i^{\mu} = \frac{d_i(\mu)}{\sqrt{D_{k-1}(\mu)D_{k}(\mu)}}. $$
Using the same notation as in the proof of Proposition \ref{propositionloclipbeta} we have that 
\begin{align*}
|a_i(\mu)- a_i(\nu)| & \leq \left| \frac{d_i(\mu)}{\sqrt{x_{k-1}x_k}}- \frac{d_i(\nu)}{\sqrt{y_{k-1}y_k}}  \right| \\ &  \leq \frac{1}{\sqrt{x_{k-1}x_k}}|d_i(\mu)-d_i(\nu)|+d_i(\nu)\left| \frac{1}{\sqrt{x_{k-1}x_k}}- \frac{1}{\sqrt{y_{k-1}y_k}}\right|.
\end{align*}

As before $\frac{1}{\sqrt{x_{k-1}x_k}}\leq \frac{1}{s_k}$, while $|d_i(\mu)-d_i(\nu)|\leq 2\mathrm{Kol}(\mu, \nu)\exp\{g k^2\}$ for some $g>0$ dependent on $\mu$ and $\nu$ only. To bound the second term on the right-hand side of the above inequality note that $d_i(\nu)\leq \exp\{g k^2\}$ and that 
\begin{align*}
\frac{1}{\sqrt{x_{k-1}x_k}}- \frac{1}{\sqrt{y_{k-1}y_k}} & = (x_{k-1}x_ky_{k-1}y_k)^{-\frac{1}{2}}|\sqrt{x_{k-1}x_k}-\sqrt{y_{k-1}y_k}| \\ &\leq \frac{1}{s_k^3} \exp\{g k^2\}\mathrm{Kol}(\mu, \nu),
\end{align*}

where the last inequality is a consequence of (\ref{ineq}). The result follows. 
\end{proof}

\begin{corollary}
\label{coralpha}
Let $\mu, \nu, s_k$ be as in Proposition \ref{propositionloclipbeta}.  Then 
$$|\alpha_k^\mu- \alpha_k^\nu|\leq \frac{\mathrm{Kol}(\mu, \nu) \exp\{g k^2\}}{s_k^{3}}. $$
\end{corollary}

\begin{proof}
Recall that 
$$\alpha_k^\mu = \int x p_k^2(x) d \mu(x) = \sum_{i, j=1}^k a_i^\mu a_j^\mu m_{i+j+1}(\mu).$$
As mentioned above, the quantities $a_i^\mu, a_i^\nu$, and $m_{i}(\mu), n_i(\nu)$ are of size $O(\exp\{gk^2\})$.
Putting this together with Proposition \ref{proplocalLipscoefficients} and Lemma \ref{lemmamomentkolmogorov}  we get that
 
$$|a_i^\mu a_j^\mu m_{i+j-1}(\mu)- a_i^\nu a_j^\nu m_{i+j-1}(\nu)| \leq \frac{\exp\{gk^2\}}{s_k^3}.$$
By adding over $i, j$ and modifying $g$ the result follows. 
\end{proof}

In order to prove Theorem \ref{Thmsqrtiter} and Proposition \ref{cormain} we need one final lemma, which states that with overwhelming probability, the random measure $\mu_n^u$ is close in Kolmogorov distance to $\mu_n$. 

\begin{lemma}
\label{lemmaKolmuu}
For $n$ large enough we have that
$$\mathbf{P}[\mathrm{Kol}(\mu_n^u, \mu_n)\geq n^{-\frac{1}{4}}] \leq \exp\{-n^{\frac{1}{4}}/8\}. $$
\end{lemma}

\begin{proof}
We must show that
$$\left| \sum_{i=1}^k u_i^2- \frac{k}{n}\right| \leq n^{-\frac{1}{4}}$$
for all $1 \le k \le n$ with probability at least $1 - \exp\{-n^{1/4}/8\}$. 

Fix $1\leq k \leq n$. As in Section 3.3 start by considering $X_1, \dots, X_k$ independent centered Gaussian random variables of variance $\frac{1}{n}$ and let $Z_k = \sum_{i=1}^k X_i^2$. Then by Lemma \ref{lemmachitail} we have that 
$$\mathbf{P}\left[Z_k \geq \frac{k}{n}+ n^{-\frac{1}{4}}\right] \leq e^{-t_1} \hspace{.3cm} \text{and} \hspace{.3cm} \mathbf{P}\left[Z_k \leq \frac{k}{n}-n^{-\frac{1}{4}}\right]\leq e^{-t_2},$$
where $t_1$ and $t_2$ are the solutions to
\begin{equation}
\label{eqt's}
n^{-\frac{1}{4}} = \frac{2\sqrt{kt_1}}{n} \hspace{.3cm} \text{and} \hspace{.3cm} n^{-\frac{1}{4}}= \frac{2\sqrt{kt_2}+2t_2}{n},
\end{equation}
respectively. Since $k\leq n$ it is clear from (\ref{eqt's}) that $\min\{t_1, t_2\} \geq \frac{n^{\frac{1}{4}}}{4}$. This implies that 
$$\mathbf{P}\left[\left|Z_k-\frac{k}{n}\right|\geq n^{-\frac{1}{4}}\right]\leq \exp\{-n^{\frac{1}{4}}/4\}. $$
Now, letting $k$ run from 1 to $n$, a union bound yields that
$$\mathbf{P}\left[\max_{1\leq k\leq n} \left|Z_k-\frac{k}{n}\right|>n^{-\frac{1}{4}} \right] \leq n \exp\{-n^{\frac{1}{4}}/4\} \leq \frac{1}{2}  \exp \{-n^{\frac{1}{4}}/8\},$$
where the last equality holds for $n$ large enough. Now, as in the proof of Proposition (\ref{propmainsphere}) we can show by a standard coupling argument that if we take  $u_i = X_i / \sqrt{Z_n} $, we will have that 
$$\mathbf{P}\left[\max_{1\leq k \leq n}\left|Z_k-\sum_{i=1}^k u_i^2\right|\right] \leq \frac{1}{2} \exp\{-n^{\frac{1}{4}}/8\}   $$
and the result follows. 
\end{proof}

\begin{proof}[Proof of Theorem \ref{Thmsqrtiter}] From Lemma \ref{lemmaKolmuu}, for  $n$ large enough,  we  have that $\mathrm{Kol}(\mu^u, \mu_n) \leq n^{-\frac{1}{4}}$ with overwhelming probability. By the assumption $\mathrm{Kol}(\mu_n, \mu) = n^{-c}$ we then have that $\mathrm{Kol}(\mu^u, \mu) \leq n^{-c'}$ also with overwhelming probability for $c' = \min\{1/4, c\}$. Hence, under the event $\{\mathrm{Kol}(\mu^u, \mu) \leq n^{-c'}\}$ we can apply Proposition \ref{propositionloclipbeta} and Corollary \ref{coralpha} and use the fact that the Jacobi matrices are tridiagonal to obtain that
$$\norm{J_{k_n}(u)-J_{k_n}(\mu)}_{{}}\leq \frac{6 C \exp\{d'k^2\}}{n^{c'} \min\{s_k^2, s_k^3\}} .   $$ 

Since $\mu$ has an absolutely continuous part we know from Proposition \ref{propositionAbcContSparse} and Corollary   \ref{corkolconv} that $|\log(\gamma_k^\mu)| = O(k)$. Hence, from (\ref{eqothgamma}) we get $|\log s_k | = O(k^2)$, which makes it clear that there exists $d> 0$ and a sequence $k_n \le d \sqrt{\log n}$ satisfying the theorem statement.
\end{proof}

\begin{proof}[Proof of Proposition \ref{cormain}] As mentioned in Section 2, this proposition is a direct consequence of Theorem \ref{Thmsqrtiter} and Lemma \ref{lemmaWielandt-Hoffman}. 
\end{proof}

\begin{remark} \label{rem:lognhard}
Observe that the above proofs repeatedly use the fact that moments are Lipschitz quantities on $\mathcal{P}_C$ and that the Jacobi coefficients are an explicit function of the moments. However, going from moments to Jacobi coefficients is an expensive process that we pay for by getting a rate of $O(\sqrt{\log n})$ instead of $\Theta (\log n)$. At first glance, it may seem that the results in Section 3.2 may be used in a similar fashion to obtain a better rate; however, even if we have strong concentration results for the Jacobi coefficients of the random measures $\mu_n^u$, it is a difficult task to control the location of the medians (or means) of $\alpha_j(u)$ and $\beta_j(u)$ and hence it is hard to show that these quantities converge at a good enough rate to the Jacobi coefficients of $\mu$.
\end{remark}

\section{Concluding remarks} \label{sec:conclusion}

Several directions can be pursued to expand the results presented throughout this paper. Currently, we have only analyzed the Lanczos algorithm in its prototypical form, but have not analyzed the more sophisticated variants that are used in practice.  Obtaining similar concentration results and negative results for these modifications, and more generally for Krylov subspace methods, would  be of great interest. 

The Lanczos algorithm is used in practice for non-Hermitian matrices and even nonnormal matrices, despite these cases being far less understood.  In this incarnation, the algorithm is referred to as the Arnoldi algorithm. Extending the results of this paper to the Arnoldi algorithm is a natural direction to pursue. As mentioned in Remark \ref{pseudospectrum}, it is easy to extend Theorem \ref{thmhessconc} to the non-Hermitian setting, but no longer so easy to prove concentration of the Ritz values or to say anything about their location.

The concentration guarantees of the output of the Lanczos algorithm in the present work hold only for up to $C\log(n)$ iterations, where $C$ is a function of the equidistribution parameters of the spectrum of the input matrix. However, we do not know what is the optimal function of $n$ for which a result of this sort holds. Proving a general concentration statement that holds for $\Omega(\log(n))$ iterations would be interesting and require essentially different ideas.  A particular case, which might have a simpler solution but is of great interest in applications, is to show concentration of the eigenvectors that correspond to outlying eigenvalues. A possible approach to the latter problem is to take classical guarantees for approximation of outlying eigenvectors, such as \cite{saad1980rates}, and convert them into probabilistic statements like the one given in Theorem \ref{thmhessconc}.

In the same direction,  another interesting task is to extend Theorem \ref{Thmsqrtiter}  to hold in the regime $k = \Omega(\log n)$ instead of the current setting of $O(\sqrt{\log n})$.  Some key difficulties are discussed in Remark \ref{rem:lognhard}.

Another direction is to translate the concentration of the Jacobi coefficients (Theorem \ref{thmhessconc}) into some quantitative statement about the quality of the approximate spectral density obtained from the first few Jacobi coefficients.  This would require analyzing the conditioning of the Haydock method or Pad\'{e} approximation mentioned in the introduction. The Jacobi coefficients are also used for estimating matrix functionals via quadrature; see \cite{gautschi2004orthogonal} for a comprehensive account.  It would be interesting to understand the implications of concentration in this setting as well.


Finally, we pose a question that is of independent mathematical interest, but which would also allow us to remove the square root on the $\log{n}$ in the statement of Theorem \ref{Thmsqrtiter} and Proposition \ref{cormain}.

\begin{question}
\label{quest1}
For $C > 0$ let $\mathcal{P}_C$ be the set of probability measures with support contained in the interval $[-C, C]$.  Is there a natural metric on $\mathcal{P}_C$ inducing a topology for which the set of atomic measures is a dense subset of $\mathcal{P}_C$, and the Jacobi coefficients $$\alpha_j : \mathcal{P}_C\to \mathbb{R} \hspace{.3cm} \text{and} \hspace{.3cm} \beta_j : \mathcal{P}_C \to \mathbb{R}$$  have local Lipschitz constant at most exponential in $j$?
\end{question}

\section*{Acknowledgments}
We would like to thank Nikhil Srivastava for many helpful discussions, for introducing us to this topic, and for providing valuable comments on a draft of this paper.  We would also like to thank Stanis\l{}aw Szarek, Lin Lin,  and Max Simchowitz for helpful discussions.  Finally, we would like to thank the anonymous referees for comments on an earlier version of this paper which helped improve the presentation.
\bibliographystyle{alpha}
\bibliography{references}

\newcommand{\etalchar}[1]{$^{#1}$}
\begin{thebibliography}{VDHVDV01}

\bibitem[AS17]{aubrun2017alice}
Guillaume Aubrun and Stanis{\l}aw~J. Szarek.
\newblock {\em Alice and Bob Meet Banach: The Interface of Asymptotic Geometric
  Analysis and Quantum Information Theory}, volume 223.
\newblock American Mathematical Society, 2017.

\bibitem[Bec00]{beckermann2000note}
Bernhard Beckermann.
\newblock A note on the convergence of {R}itz values for sequences of matrices.
\newblock {\em Publication ANO}, 408, 2000.

\bibitem[BSS10]{bellalij2010further}
Mohammed Bellalij, Yousef Saad, and Hassane Sadok.
\newblock Further analysis of the {A}rnoldi process for eigenvalue problems.
\newblock {\em SIAM J. Numer. Anal.}, 48(2):393--407, 2010.

\bibitem[CRS94]{calvetti1994implicitly}
Daniela Calvetti, Lothar Reichel, and Danny~Chris Sorensen.
\newblock An implicitly restarted {L}anczos method for large symmetric
  eigenvalue problems.
\newblock {\em Electron. Trans. Numer. Anal.}, 2(1):21, 1994.

\bibitem[Dei99]{deift1999orthogonal}
Percy Deift.
\newblock {\em Orthogonal polynomials and random matrices: a
  {R}iemann-{H}ilbert approach}, volume~3.
\newblock American Mathematical Soc., 1999.

\bibitem[DK69]{davis1969some}
Chandler Davis and William~M Kahan.
\newblock Some new bounds on perturbation of subspaces.
\newblock {\em Bull. Amer. Math. Soc.}, 75(4):863--868, 1969.

\bibitem[Gau68]{gautschi1968construction}
Walter Gautschi.
\newblock Construction of gauss-christoffel quadrature formulas.
\newblock {\em Math. Comp.}, 22(102):251--270, 1968.

\bibitem[Gau04]{gautschi2004orthogonal}
Walter Gautschi.
\newblock {\em Orthogonal polynomials: computation and approximation}.
\newblock Oxford University Press on Demand, 2004.

\bibitem[God17]{godsil2017algebraic}
Chris Godsil.
\newblock {\em Algebraic combinatorics}.
\newblock Routledge, 2017.

\bibitem[GU77]{golub1977block}
Gene~Howard Golub and Richard Underwood.
\newblock The block {L}anczos method for computing eigenvalues.
\newblock In {\em Mathematical software}, pages 361--377. Elsevier, 1977.

\bibitem[Hay80]{haydock1980recursive}
Roger Haydock.
\newblock The recursive solution of the schr{\"o}dinger equation.
\newblock {\em Comput. Phys. Commun.}, 20(1):11--16, 1980.

\bibitem[HJ12]{horn2012matrix}
Roger~A Horn and Charles~R Johnson.
\newblock {\em Matrix analysis}.
\newblock Cambridge University Press, 2012.

\bibitem[Kan66]{kaniel1966estimates}
Shmuel Kaniel.
\newblock Estimates for some computational techniques in linear algebra.
\newblock {\em Math. Comp.}, 20(95):369--378, 1966.

\bibitem[Kui00]{kuijlaars2000eigenvalues}
Arno B.~J. Kuijlaars.
\newblock Which eigenvalues are found by the {L}anczos method?
\newblock {\em SIAM J. Matrix Anal. Appl.}, 22(1):306--321, 2000.

\bibitem[Kui06]{kuijlaars2006convergence}
Arno B.~J. Kuijlaars.
\newblock Convergence analysis of {K}rylov subspace iterations with methods
  from potential theory.
\newblock {\em SIAM Rev.}, 48(1):3--40, 2006.

\bibitem[KW94]{kuczynski1994probabilistic}
Jacek Kuczy{\'n}ski and Henryk Wo{\'z}niakowski.
\newblock Probabilistic bounds on the extremal eigenvalues and condition number
  by the {L}anczos algorithm.
\newblock {\em SIAM J. Matrix Anal. Appl.}, 15(2):672--691, 1994.

\bibitem[LM00]{laurent2000adaptive}
Beatrice Laurent and Pascal Massart.
\newblock Adaptive estimation of a quadratic functional by model selection.
\newblock {\em Ann. Statist.}, pages 1302--1338, 2000.

\bibitem[LSY16]{lin2016approximating}
Lin Lin, Yousef Saad, and Chao Yang.
\newblock Approximating spectral densities of large matrices.
\newblock {\em SIAM Rev.}, 58(1):34--65, 2016.

\bibitem[LXV{\etalchar{+}}16]{li2016thick}
Ruipeng Li, Yuanzhe Xi, Eugene Vecharynski, Chao Yang, and Yousef Saad.
\newblock A {T}hick-{R}estart {L}anczos algorithm with polynomial filtering for
  {H}ermitian eigenvalue problems.
\newblock {\em SIAM J. Sci. Comput.}, 38(4):A2512--A2534, 2016.

\bibitem[Mul59]{muller1959note}
Mervin~E Muller.
\newblock A note on a method for generating points uniformly on n-dimensional
  spheres.
\newblock {\em Commun. ACM}, 2(4):19--20, 1959.

\bibitem[Pai71]{paige1971computation}
Christopher~Conway Paige.
\newblock {\em The computation of eigenvalues and eigenvectors of very large
  sparse matrices}.
\newblock PhD thesis, University of London, 1971.

\bibitem[Saa80]{saad1980rates}
Yousef Saad.
\newblock On the rates of convergence of the {Lanczos} and the {block-Lanczos}
  methods.
\newblock {\em SIAM J. Numer. Anal.}, 17(5):687--706, 1980.

\bibitem[Saa11]{saad2011numerical}
Yousef Saad.
\newblock {\em Numerical methods for large eigenvalue problems: revised
  edition}, volume~66.
\newblock SIAM, 2011.

\bibitem[SdJL{\etalchar{+}}18]{shao2018structure}
Meiyue Shao, Felipe~H da~Jornada, Lin Lin, Chao Yang, Jack Deslippe, and
  Steven~G Louie.
\newblock A structure preserving lanczos algorithm for computing the optical
  absorption spectrum.
\newblock {\em SIAM J. Matrix Anal. Appl.}, 39(2):683--711, 2018.

\bibitem[SEAR18]{simchowitz2018tight}
Max Simchowitz, Ahmed El~Alaoui, and Benjamin Recht.
\newblock Tight query complexity lower bounds for {PCA} via finite sample
  deformed {Wigner} law.
\newblock In {\em Proceedings of the 50th Annual ACM SIGACT Symposium on Theory
  of Computing}, pages 1249--1259. ACM, 2018.

\bibitem[Sze39]{szego1939orthogonal}
Gabor Szeg\H{o}.
\newblock {\em Orthogonal polynomials}, volume~23.
\newblock American Mathematical Society, 1939.

\bibitem[Sze77]{Szego}
Gabor Szeg\H{o}.
\newblock Hankel forms.
\newblock {\em Amer. Math. Soc. Transl.}, 108, 1977.

\bibitem[TBI97]{trefethen1997numerical}
Lloyd~N. Trefethen and David Bau~III.
\newblock {\em Numerical linear algebra}, volume~50.
\newblock SIAM, 1997.

\bibitem[VA06]{van2006pade}
Walter Van~Assche.
\newblock Pad{\'e} and hermite-pad{\'e} approximation and orthogonality.
\newblock {\em Surv. Approx. Theory}, 2:61--91, 2006.

\bibitem[VDHVDV01]{van2001computing}
Jos L.~M. Van~Dorsselaer, Michiel~E. Hochstenbach, and Henk~A. Van Der~Vorst.
\newblock Computing probabilistic bounds for extreme eigenvalues of symmetric
  matrices with the {L}anczos method.
\newblock {\em SIAM J. Matrix Anal. Appl.}, 22(3):837--852, 2001.

\bibitem[Ver09]{vershynin2009role}
Roman Vershynin.
\newblock On the role of sparsity in compressed sensing and random matrix
  theory.
\newblock In {\em 2009 3rd IEEE International Workshop on Computational
  Advances in Multi-Sensor Adaptive Processing (CAMSAP)}, pages 189--192. IEEE,
  2009.

\bibitem[Ver18]{vershynin2018high}
Roman Vershynin.
\newblock {\em High-dimensional probability: An introduction with applications
  in data science}, volume~47.
\newblock Cambridge University Press, 2018.

\bibitem[YGL18]{yuan2018superlinear}
Qiaochu Yuan, Ming Gu, and Bo~Li.
\newblock Superlinear convergence of randomized block {L}anczos algorithm.
\newblock In {\em 2018 IEEE International Conference on Data Mining}, pages
  1404--1409. IEEE, 2018.

\end{thebibliography}
\end{document}